\documentclass[11pt]{amsart}

\usepackage{amsmath,amsthm, amscd, amssymb, amsfonts}
\usepackage[all]{xy}
\usepackage[dvips, dvipsnames, usenames]{color}

\newtheorem{theorem}{Theorem}[section]

\newtheorem{lemma}[theorem]{Lemma}
\newtheorem{corollary}[theorem]{Corollary}

\newtheorem{teo intro}{Theorem}

\newtheorem{proposition}[theorem]{Proposition}

\theoremstyle{definition}
\newtheorem{definition}[theorem]{Definition}

\theoremstyle{remark}
\newtheorem{remark}[theorem]{Remark}

\def\yd{^{H}_{H}\mathcal{YD}}

\newcommand{\ydzt}{^{\ku\zt}_{\ku\zt}\mathcal{YD}}

\def\zt{\Z^{\theta}}

\newcommand{\g}{\mathfrak g}
\newcommand\id{\operatorname{id}}

\newcommand\ord{\operatorname{ord}}
\newcommand\Hom{\operatorname{Hom}}

\newcommand\hgt{\operatorname{ht}}

\newcommand\Aut{\operatorname{Aut}}

\newcommand\cop{\operatorname{cop}}
\newcommand\bop{\operatorname{bop}}

\newcommand\gr{\operatorname{gr}}
\newcommand\co{\operatorname{co}}

\newcommand\ad{\operatorname{ad}}
\newcommand\Spec{\operatorname{Spec}}
\newcommand\GKdim{\operatorname{GKdim}}

\newcommand\re{\operatorname{re}}

\newcommand{\ku}{ \mathbf{k}}

\def\ot{\otimes}

\def\C{\mathbb{C}}
\def\Z{\mathbb{Z}}
\def\N{\mathbb{N}}

\def\cB{\mathcal{B}}
\def\toba{\cB}
\def\dpn{\widetilde{\mathcal{B}}}
\def\cO{\mathcal{O}}

\def\cC{\mathcal{C}}
\def\cI{\mathcal{I}}
\def\cJ{\mathcal{J}}
\def\cK{\mathcal{K}}
\def\cX{\mathcal{X}}
\def\cR{\mathcal{R}}
\def\cQ{\mathcal{Q}}

\def\cG{\mathcal{G}}
\def\cH{\mathcal{H}}
\def\cW{\mathcal{W}}

\def\cP{\mathcal{P}}
\def\u{\mathfrak{u}}

\def\Ef{\mathfrak{E}}
\def\Rf{\mathfrak{R}}

\newcommand\Sb{\mathbb S}

\def\cU{\mathcal{U}}

\newcommand{\bm}{\mathbf{m} }
\newcommand{\bM}{\mathbf{M} }

\newcommand\I{\mathbb I}

\newcommand{\Eb}{\underline E}
\newcommand{\Fb}{\underline F}
\newcommand{\Kb}{\underline K}
\newcommand{\Lb}{\underline L}

\newcommand{\cb}{\underline c}

\newcommand{\E}{\mathbf E}
\newcommand{\F}{\mathbf F}
\newcommand{\Et}{\mathtt E}
\newcommand{\Ft}{\mathtt F}

\newcommand{\ab}{\mathbf{a}}
\newcommand{\bb}{\mathbf{b}}
\newcommand{\ub}{\mathbf{u}}
\newcommand{\vb}{\mathbf{v}}

\newcommand{\schi}{{\rho_i(V)}}

\newcommand{\eps}{\varepsilon}
\newcommand{\un}{_{(1)}}
\newcommand{\dos}{_{(2)}}

\newcommand{\cero}{_{(0)}}
\newcommand{\mun}{_{(-1)}}

\newcommand{\tri}{\triangleright}

\newcommand{\vi}{${\mathsf {(i)}\;}$}
\newcommand{\vii}{${\mathsf {(ii)}\;}$}
\newcommand{\viii}{${\mathsf {(iii)}\;}$}

\def\pf{\begin{proof}}
\def\epf{\end{proof}}

\begin{document}


\title[Distinguished Pre-Nichols algebras]{Distinguished Pre-Nichols algebras}
\author[Iv\'an Angiono]{Iv\'an Angiono}

\address{FaMAF-CIEM (CONICET), Universidad Nacional de C\'ordoba,
Medina A\-llen\-de s/n, Ciudad Universitaria, 5000 C\' ordoba, Rep\'
ublica Argentina.} \email{angiono@famaf.unc.edu.ar}

\thanks{The work was partially supported by CONICET,
FONCyT-ANPCyT, Secyt (UNC)}

\begin{abstract}
We formally define and study the distinguished pre-Nichols algebra $\dpn(V)$
of a braided  vector space of diagonal type $V$ with finite-dimensional Nichols algebra $\cB(V)$.
The algebra $\dpn(V)$ is presented by fewer relations than $\cB(V)$, so it is
intermediate between the tensor algebra $T(V)$ and $\cB(V)$. Prominent examples of distinguished pre-Nichols algebras
are the positive parts of quantized enveloping (super)algebras and their multiparametric versions.
We prove that these algebras give rise to new examples of
Noetherian pointed Hopf algebras of finite Gelfand-Kirillov dimension.
We investigate the kernel (in the sense of Hopf algebras) of the projection from $\dpn(V)$ to $\cB(V)$,
generalizing results of De Concini and Procesi on quantum groups at roots of unity.
\end{abstract}

\maketitle

\section{Introduction}

\subsection{\hspace{0.1pt} }
Let $\g$ be a finite-dimensional semisimple Lie algebra over $\C$.
The quantized enveloping algebra $U_q(\g)$ was introduced in \cite{Dr-quantum,J}
by deforming the relations of the enveloping algebra $U(\g)$.
If the parameter $q$ is not a root of unity, then $U_q(\g)$ (defined as in \cite{DP-quantum}) has similar properties to $U(\g)$.
But if $q$ is a root of unity, then $U_q(\g)$ differs significantly from $U(\g)$,
For instance
there is a Hopf pairing between $U_q^{+}(\g)$ (the subalgebra generated by the $E_i$'s)
and $U_q^{-}(\g)$ (the subalgebra generated by the $F_i$'s)  but is degenerate
and its radical is the ideal generated by powers of roots vectors $E_\alpha^N$,
respectively $F_\alpha^N$,  $\alpha\in\Delta_+$. Here $E_\alpha$, $F_\alpha$,  $\alpha\in\Delta_+$,
are obtained from $E_i$, $F_i$ by applying Lusztig isomorphisms $T_i\in \Aut U_q(\g)$ \cite[Section 9]{DP-quantum}.
The quotient of $U_q(\g)$ by the ideal generated by $E_\alpha^N$, $K_\alpha^N$, $F_\alpha^N$ is a finite-dimensional Hopf algebra
known as the \emph{small quantum group} or \emph{Frobenius-Lusztig kernel} $\u_q(\g)$.
Now the induced Hopf pairing between $\u_q^{+}(\g)$ and $\u_q^{-}(\g)$ is non-degenerate, so that $\u_q^{+}(\g)$ is a Nichols algebra
\cite{AS Pointed HA}.
The kernel of the natural projection $U_q(\g) \to \u_q(\g)$ (in the  Hopf algebra sense) is the central Hopf subalgebra $Z_0$ generated
by $E_\alpha^N$, $K_\alpha^N$, $F_\alpha^N$; $Z_0$ is the algebra of functions of a Poisson group \cite[Chapter 19]{DP-quantum}
that plays roles in the representation theory of $U_q(\g)$ \cite{DP-quantum, DPRR} and in the classification of finite-dimensional pointed Hopf algebras
\cite{AS Class}.

\subsection{\hspace{0.1pt} }
A \emph{braided vector space} is a pair $(V,c)$, where $V$ is a vector space and $c\in\Aut(V\ot V)$ satisfies
$(c\ot \id )(\id \ot c)(c\ot \id )=(\id \ot c)(c\ot \id )(\id \ot c)$. The Nichols algebra  of a braided vector space $(V,c)$
is a graded braided Hopf algebra $\cB(V) = T(V)/\cJ(V) = \oplus_{n\ge 0}\cB^n(V)$ with the property
that all primitive elements are in degree 1.
The study of Nichols algebras is crucial in the classification program of Hopf algebras \cite{AS Pointed HA}.
If there are a basis $(v_i)_{1\le i\le\theta}$ of $V$ and  a matrix $(q_{ij})\in (\C^{\times})^{\theta\times \theta}$
such that $c(v_i\ot v_j)=q_{ij} \, v_j\ot v_i$ for all $1\le i,j\le\theta$, then  $(V,c)$ is of \emph{diagonal type}.
Finite-dimensional Nichols algebras of (braided vector spaces of) diagonal type are classified in \cite{H-classif RS};
the defining relations of them are described in \cite{A-convex,A-presentation}.

A pre-Nichols algebra\footnote{This terminology is due to Akira Masuoka.} of a braided vector space $(V,c)$ is any graded braided Hopf algebra intermediate between $T(V)$ and $\cB(V)$,
that is any braided Hopf algebra of the form $T(V)/\cI$ where $\cI \subseteq \cJ(V)$ is a homogeneous Hopf ideal.
Pre-Nichols algebras form a partially ordered set with order given by projection,
with $T(V)$ the minimal and $\toba(V)$ the maximal points.
Pre-Nichols algebras appear naturally in the computation of the deformations or liftings of \cite{Masuoka, AAGMV}.
Let $V$ be a braided  vector space of diagonal type  with finite-dimensional  $\cB(V)$.
In this paper we define and investigate the \emph{distinguished} pre-Nichols algebra $\dpn(V)$.
Actually they are already present without name in \cite[Proposition 3.3]{A-presentation}.
The distinguished pre-Nichols algebra $\dpn(V)$ can be realized in the category of Yetter-Drinfeld modules
over $\ku\zt$, where $\theta = \dim V$; let $U(V)$ be the quantum double of the bosonization
$\dpn(V)\# \ku\zt$. Similarly, let $\u(V)$ be the quantum double of the bosonization
$\cB(V)\# \ku\zt$.
Here are some properties of $\dpn(V)$, justifying the adjective distinguished (to our understanding):	

\subsubsection{} Let $i\in\I_{\theta}$ and $\rho_i: V \mapsto V'$ the corresponding reflection in the Weyl groupoid of $(V,c)$.
The Lusztig isomorphism  $T_i:\u(V) \to \u(V')$ can be lifted to a Lusztig isomorphism $T_i: U(V)\to U(V')$ \cite[Proposition 3.26]{A-presentation}.
We conjecture that $\dpn(V)$ is minimal among the pre-Nichols algebras admitting Lusztig isomorphisms.

\subsubsection{} If $\toba (V) = \u_q^{+}(\g)$, then $\dpn(V) = U_q^{+}(\g)$.

\subsubsection{} By a general result of \cite{Kh}, every pre-Nichols algebra has a restricted (that is, with heights) PBW basis.
We prove the existence of a PBW basis with the same generators (root vectors), but different heights,
as $\toba(V)$, obtained by applying several $T_i$'s, see Theorem \ref{thm: PBW bases preNichols}.

\subsubsection{} The pointed Hopf algebra $U(V)$ is Noetherian of finite Gelfand-Kiri\-llov dimension, see Theorems \ref{thm:noetheriano} and \ref{thm:GKdim finita}.
We conjecture\footnote{jointly with N. A.} that $\dpn(V)$ is minimal among the pre-Nichols algebras with these properties.
Observe that $\dpn(V)$ is not a domain in general.

\subsubsection{} The powers of root vectors that are non-zero in the $\dpn(V)$ but zero in $\toba(V)$
generate a subalgebra $Z^+(V)$ on $\dpn(V)$, that coincides with the intersection of the kernels of the
skew-derivations associated to the coproduct of $U^+(V)$, see Theorem \ref{thm:Z+ es subalg hopf U+}. Correspondingly there is a normal Hopf subalgebra $Z(V)$ of $U(V)$,
see Theorem \ref{thm:Z es subalg hopf U};
$U(V)$ is a finite free $Z(V)$-module.

\subsection{\hspace{0.1pt} }The organization of the article is the following. In Section \ref{section:preliminares} we recall notions and
basic properties of generalized root systems and Lusztig isomorphisms of quantum doubles of Nichols algebras.

In Section \ref{section:Distinguished pre-Nichols algebras} we study the distinguished pre-Nichols algebras $\dpn(V)$.
First we prove the existence of Lusztig isomorphisms for the quantum double $U(V)$ and PBW bases
whose PBW generators are obtained by applying Lusztig isomorphisms. Then we obtain an algebra filtration on $\dpn(V)$ and $U(V)$ such that the associated graded algebra is a quantum polynomial algebra, so $\dpn(V)$ and $U(V)$ are Noetherian of finite Gelfand-Kirillov dimension.

In Section \ref{section:subalgebra Z} we study the subalgebra $Z(V)$ of $U(V)$ generated by powers of root vectors.
First we show that each element of $Z(V)$ commutes up to scalars with each homogeneous element of $U(V)$, so $U(V)$ is a free $Z(V)$-module.
Next we obtain a formula relating the coproduct on $\dpn(V)$ with Lusztig isomorphisms, close to \cite[Theorem 4.2]{HS-coideal subalg}.
We present a recursive formula for the coproduct of powers of root vectors in $\dpn(V)$, with a view towards
the computation of the liftings of $\toba(V)$, as proposed in \cite{AAGMV}.

Finally, we compute the coproduct of $Z(V)$ for braidings of super type $A$ and of type $\mathfrak{br}(2;5)$ in Section \ref{section:examples}.

\subsection{\hspace{0.1pt} }
Here are some questions on, and potential applications of, distinguished pre-Nichols algebras that support our interest
on them.

\subsubsection{} If the diagonal braiding $(V,c)$ is given by a symmetric matrix, then the algebra $Z(V)$ is commutative.
Let $G(V)$ be the algebraic group  $\Spec Z(V)$.  There exists a correspondence
between (almost all) braidings of diagonal type whose Nichols algebras are finite-dimensional,
and finite-dimensional contragredient  Lie superalgebras in positive characteristic \cite{AA-WGCLSandNA}.
The group $G(V)$ should be related with the Lie superalgebra corresponding to $(V,c)$;
it should allow a geometric approach to the representation theory of $U(V)$ as in \cite{DP-quantum}.

\subsubsection{} The Hopf algebra $U(V)$ is the `quantum' analogue of the enveloping Lie superalgebra in the previous point
while $\u(V)$ is the one of its restricted enveloping superalgebra.
Therefore, as for Lie algebras in positive characteristic,
we can start by studying the representation theory of the algebra $U(V)$ and pass to $\u(V)$.

\subsection*{Acknowledgements}
I thank Nicol\'as Andruskiewitsch for interesting and guiding discussions, and several comments which help to improve this work.

\section{Preliminaries}\label{section:preliminares}

\subsection{Conventions}

We work over an algebraically closed field $\ku$ of characteristic zero. Tensor products, algebras, coalgebras and Hopf algebras are taken over $\ku$. For each $\theta \in\N$, let $\I_\theta=\{1,2,\dots,\theta\}$; when $\theta$ is clear from the context, we simply set $\I=\I_\theta$.

Given a matrix $(q_{ij})_{1\le i,j\le\theta}\in\ku^{\theta\times\theta}$, let $\widetilde{q_{ij}}=q_{ij}q_{ji}$, $i\neq j$.

We consider the $q$-polynomial numbers in $\Z[\mathbf{q}]$, $ n\in \N$, $0 \leq i \leq n$,
\begin{align*}
(n)_\mathbf{q} &=\sum_{j=0}^{n-1}\mathbf{q}^{j}, & (n)_\mathbf{q}!&=\prod_{j=1}^{n} (j)_\mathbf{q}, &
\binom{n}{i}_\mathbf{q} & =\frac{(n)_\mathbf{q}!}{(n-i)_\mathbf{q}!(i)_\mathbf{q}!}.
\end{align*}
$(n)_q$, $(n)_q!$, $\binom{n}{i}_q$ are the evaluations of the polynomials in $\mathbf{q}=q\in\ku$.

Given $\ab=(a_1,\dots,a_\theta)\in\N_0^\theta$, we use the notation: $\mathbf{t}^\ab= t_1^{a_1}\cdots t_\theta^{a_\theta}$.
The \emph{height} of $\ab$ is $\hgt(\ab)=\sum_{j=1}^\theta a_j$.
If $W=\oplus_{\ab\in\N_0^\theta} W_\ab$ is an $\N_0^\theta$-graded vector space, then its \emph{Hilbert series} is $ \mathfrak{H}_W = \sum_{\ab\in\N_0^\theta} \dim W_\ab \, \mathbf{t}^\ab \in\N_0[[t_1,\dots,t_\theta]]$.

We denote by $\{\alpha_i\}_{1\le i\le\theta}$ the canonical basis on $\Z^\theta$. Given a bicharacter $\chi:\zt\times\zt\to\ku^\times$ and $w\in\Aut(\zt)$, $w^*\chi$ denotes the bicharacter
\begin{align*}
w^*\chi(\beta,\gamma)&=\chi(w^{-1}(\beta),w^{-1}(\gamma)), & \beta,\gamma & \in\zt.
\end{align*}

\medbreak

Let $H$ be a Hopf algebra with bijective antipode $\mathcal{S}$.
We use the Sweedler notation for the comultiplication, $\Delta(h)=h\un\ot h\dos$, $h\in H$, and for left $H$-comodules $X$, $\lambda(x)=x\mun\ot x\cero\in X\ot H$, $x\in X$.

We denote by $\yd$ the category of left Yetter-Drinfeld modules over $H$. Recall that this is a braided tensor category, where the braiding for $M,N \in\yd$ is $c=c_{M,N}:M\ot N \to N\ot M$,
\begin{align*}
c(m\ot n) &= m\mun\cdot n \ot m\cero, & \quad m\in M, \, &n\in N.
\end{align*}
The inverse braiding is $c^{-1}=c_{M,N}^{-1}:N\ot M \to M\ot N$,
\begin{align*}
c^{-1}(n\ot m)&= m\cero\ot \mathcal{S}^{-1}(m\mun)\cdot n, & \quad m\in M, \, &n\in N.
\end{align*}
In particular, if $V\in\yd$, then $(V,c_{V\ot V})$  is a braided vector space.
If $(V,c)$ is of diagonal type with basis $E_1,\dots,E_\theta$ and matrix $(q_{ij})$, $\Gamma$ is a group
and $g_i\in\Gamma$, $\chi_j\in\widehat\Gamma$
are such that $\chi_j(g_i)=q_{ij}$, then $(V,c)$ can be realized as a Yetter-Drinfeld module over $H=\ku\Gamma$ by defining the action and coaction:
\begin{align}
\lambda(E_i)&=g_i\ot E_i, & g\cdot E_i&=\chi_i(g) E_i, & 1\le i\le\theta, \, & g\in\Gamma. \label{eq:YD diagonal}
\end{align}
For example take $\Gamma=\Z^\theta$, $g_i=\alpha_i$ and $\chi_j\in\widehat\zt$ such that $\chi_j(\alpha_i)=q_{ij}$ \cite{AS Pointed HA}.

Let $R$ be a Hopf algebra in $\yd$ with product $m$, coproduct $\underline{\Delta}$ and antipode $\mathcal{S}_R$. We use the following variation of Sweedler notation: $\underline{\Delta}(r)=r^{(1)}\ot r^{(2)}$, $r\in R$. $R\# H$ denotes the \emph{bosonization} of $R$ by $H$; we refer to \cite[Section 1.5]{AS Pointed HA} for the definition.

The \emph{braided commutator} and the \emph{braided adjoint action} of $x$ on $y$, $x,y\in R$, are defined as follows
\begin{align*}
[x,y]_c&=m (\id_{R\ot R}-c)(x\ot y), \\ \ad_c x(y)&=m(m\ot \mathcal{S}_R)(\id\ot c)(\underline{\Delta}\ot\id)(x\ot y).
\end{align*}
The set of primitive elements is $\cP(R)=\{r\in R: \underline{\Delta}(r)=r\ot 1+1\ot r \}$.
Notice that $\ad_c x(y)=[x,y]_c$ if $x\in\cP(R)$.

\subsection{Normal coideal subalgebras and Hopf ideals}

A \emph{right coideal subalgebra} $B$ of a Hopf algebra $H$ is a subalgebra $B$ such that $\Delta(B)\subset B\ot H$. A right coideal subalgebra $B$ is \emph{normal} if it is stable under the right adjoint action $\ad_r(h) (b)= \mathcal{S}(h_{(1)})b h_{(2)}$, $b\in B$, $h\in H$.
If $B\subset H$ is a right coideal subalgebra, then the \emph{normal right coideal subalgebra} $N(B)$\label{def:NX} is the subalgebra generated by $\{\mathcal{S}(h\un)b h\dos:\,h\in H,
 b\in B\}$.

There is a correspondence between quotient Hopf algebras and normal right coideal subalgebras under some mild conditions.

\begin{theorem}\label{thm:takeuchi-correspondance}\cite[Theorem 3.2]{T}
The maps
\begin{align*}
&B \text{ a normal right coideal subalgebra}\rightsquigarrow {\mathcal I}(B)=H B^+, \\
&I \text{ a Hopf ideal } \rightsquigarrow {\mathcal X}(I)=H^{\co H/I},
\end{align*}
restrict to mutually inverse bijective correspondences between the set of normal
right coideal subalgebras $B$ such that $H$ is right $B$-faithfully flat and the set of Hopf ideals $I$ such that $H$ is $H/I$-coflat.\qed
\end{theorem}

There is an analogous bijection for connected Hopf algebras in ${}^H_H\mathcal{YD}$.

\begin{proposition}\cite[Proposition 3.6]{AAGMV}\label{pro:masuoka}
Let $R$ be a connected Hopf algebra in
${}^H_H\mathcal{YD}$.

\smallbreak
\noindent\vi
$R$ is a free left and right module over every right coideal subalgebra, and is a cofree left and right comodule over every quotient left module coalgebra.

\smallbreak
\noindent\vii
The maps $B \mapsto R/RB^+$, $T \mapsto {}^{\co T}R$ give a bijection between the
set of right coideal subalgebras $B$ of $R$ and the set of quotient left
$R$-module coalgebras $T$ of $R$.

\smallbreak
\noindent\viii
If $B$ and $T$ correspond to each other via the bijection in \vii, then there
exists a left $T$-colinear and right $B$-linear isomorphism $T\otimes B \simeq R$.
\qed
\end{proposition}




\subsection{Weyl groupoids and generalized root systems}\label{subsection:weylgrupoide}

We recall the notion of generalized root systems using the notation in \cite{AA-WGCLSandNA}, see also \cite{CH1}.

\subsubsection{Basic data}
Given $\I=\I_\theta$, $\cX\neq\emptyset$ and $\rho: \I  \to \Sb_{\cX}$, the pair $(\cX, \rho)$ is called a \emph{basic datum}
of rank $\vert\cX\vert$ and type $\theta$ if $\rho_i^2 = \id$ for all $i\in \I$.

Let $\cQ_{\rho}$ be the quiver $\{\sigma_i^x := (x, i, \rho_i(x)): i\in \I, x\in \cX\}$ over $\cX$, with target and source
$t(\sigma_i^x) = x$, $s(\sigma_i^x) = \rho_i(x)$.
In any quotient of the free groupoid $F(\cQ_{\rho})$, we adopt the convention
\begin{align}\label{eq:conventio}
\sigma_{i_1}^x\sigma_{i_2}\cdots \sigma_{i_t} =
\sigma_{i_1}^x\sigma_{i_2}^{\rho_{i_1}(x)}\cdots \sigma_{i_t}^{\rho_{i_{t-1}} \cdots \rho_{i_1}(x)}.
\end{align}
That is, the implicit superscripts are the only possible to have compositions.

\subsubsection{Coxeter groupoids}
A \emph{Coxeter datum} is a triple $(\cX, \rho, \bM)$, where $(\cX, \rho)$ is a basic datum of type $\I$
and $\bM = (\bm^x)_{x\in \cX}$, $\bm^x = (m^x_{ij})_{i,j\in \I}$, is a family of Coxeter matrices such that
\begin{align}\label{eq:coxeter-datum}
s((\sigma^x_i\sigma_j)^{m^x_{ij}}) &= x,& i,j&\in \I,& x&\in \cX.
\end{align}
The \emph{Coxeter groupoid} $\cW(\cX, \rho, \bM)$ \cite[Definition 1]{HY} is the groupoid generated by $\cQ_{\rho}$ with relations
\begin{align}\label{eq:def-coxeter-gpd}
 (\sigma_i^x\sigma_j)^{m^x_{ij}} &= \id_x, &i, j\in \I,\,  &x\in \cX.
\end{align}
In particular, if $i = j$, then \eqref{eq:def-coxeter-gpd} means that either $\sigma_i^x$ is an involution when $ \rho_i(x) = x$, or else that $\sigma_i^x$ is the inverse arrow of
$\sigma_i^{\rho_i}(x)$ when $\rho_i(x) \neq x$.

\subsubsection{Generalized root systems}

Let $\cC=(C^x)_{x\in \cX}$ be a family of generalized Cartan matrices $C^x=(c_{ij}^x)_{i,j\in\I}$ with row invariance
\begin{align}\label{eq:condicion Cartan scheme}
c^x_{ij}&=c^{\rho_i(x)}_{ij} &  \mbox{for all }&x \in \cX, \, i,j \in \I.
\end{align}
We define $s_i^x\in GL_{\theta}(\Z)$ by
\begin{align}\label{eq:reflection-x}
s_i^x(\alpha_j)&=\alpha_j-c_{ij}^x\alpha_i, & j&\in \I,&  i &\in \I, x \in \cX.
\end{align}
Then \eqref{eq:condicion Cartan scheme} means that $s_i^x$ is the inverse of $s_i^{\rho_i(x)}$. A \emph{generalized root system} (GRS for short)
\cite[Definition 1]{HY} is a collection $\cR:= \cR(\cX, \rho, \cC, \Delta)$, where $(\cX, \rho)$
is a basic datum of type $\I$, $\cC$ is as above,
and  $\Delta = (\Delta^x)_{x\in \cX}$ is a family of subsets $\Delta^x \subset \Z^{\I}$ such that for all $x \in \cX$, $i \neq j \in \I$,
\begin{align}
\label{eq:def root system 1}
\Delta^x &= \Delta^x_+ \cup \Delta^x_-, & \Delta^x_+ &:= \Delta^x \cap \N_0^I, \,\, \Delta^x_-:=-\Delta^x_+;
\\ \label{eq:def root system 2}
\Delta^x \cap \Z \alpha_i &= \{\pm \alpha_i \};& &
\\ \label{eq:def root system 3}
s_i^x(\Delta^x)&=\Delta^{\rho_i(x)};&&
\\ \label{eq:def root system 4}
(\rho_i\rho_j)^{m_{ij}^x}(x)&=(x), &  m_{ij}^x&:=|\Delta^x \cap (\N_0\alpha_i+\N_0 \alpha_j)|.
\end{align}
We call $\Delta^x_+$, $ \Delta^x_-$ the set of \emph{positive}, respectively \emph{negative}, roots.

Let $\bM=(M^x)_{x\in\cX}$, $M^x=(m_{ij}^x)_{i,j\in\I}$.
The \emph{Weyl groupoid} of $\cR$ is $\cW = \cW(\cX, \rho, \bM)$. We can describe this groupoid by using \cite[Theorem 1]{HY}.
Let $\cG = \cX \times GL_{\theta}(\Z) \times \cX$, $\varsigma_i^x = (x, s_i^x,\rho_i(x))$, $i \in \I$, $x \in \cX$,
and $\cW' = \cW(\cX, \rho, \cC)$ the subgroupoid of $\cG$ generated by all the $\varsigma_i^x$. There exists a morphism of quivers
$\cQ_{\rho} \to \cG$, $\sigma_i^x \mapsto \varsigma_i^x$ with image $\cW'$. This induces an
isomorphism of groupoids $\cW\to \cW(\cX, \rho, \cC)$.

If $w = \sigma_{i_1}^x \cdots \sigma_{i_m}$ and $\alpha\in\zt$, then we define $w(\alpha) = s_{i_1}^x \cdots s_{i_m}(\alpha)$, so
$w(\Delta^x)= \Delta^y$, by \eqref{eq:def root system 3}. The set of \emph{real roots} at $x$ is
\begin{align}
\label{defrealroot} (\Delta^{\re})^x &= \bigcup_{y\in \cX}\{ w(\alpha_i): \ i \in \I, \ w \in \cW(y,x) \}.
\end{align}
The \emph{length} of $w\in \cW(x,\cX)$ is
$$ \ell(w)= \min \{ m\in \N_0: \ \exists i_1, \ldots, i_n \in \I \mbox{ such that }w = \sigma_{i_1}^x \cdots \sigma_{i_m} \}. $$
An expression $w = \sigma_{i_1}^x \cdots \sigma_{i_m}$ is \emph{reduced} if $m = \ell(w)$.

\begin{lemma}{\cite[Corollary 3]{HY}}  \label{Lemma:longitudHY}
Let $m \in \N$, $x, y \in \cX$, $i_1, \ldots, i_m,j \in I$, $w=\sigma_{i_1}^x \cdots \sigma_{i_m} \in \Hom(Y,X)$, where $\ell (w)=m$. Then
\begin{itemize}
  \item $\ell (w \sigma_j)=m+1$ if and only if $w(\alpha_j) \in \Delta^x_+$,
  \item $\ell (w \sigma_j)=m-1$ if and only if $w(\alpha_j) \in \Delta^x_-$.\qed
\end{itemize}
\end{lemma}

\begin{proposition}{\cite[Prop. 2.12]{CH1}}  \label{Prop:maxlongitudCH}
If $w=\sigma_{i_1}^x \cdots \sigma_{i_N}\in\cW$ is such that
$\ell(w)=N$, then all the roots $\beta_j=s_{i_1}^x\cdots
s_{i_{j-1}}(\alpha_{i_j})\in\Delta^x$ are positive and pairwise
different. In particular, if $\cR$ is finite and $w$ is an element of maximal
length, then $\Delta^x_+= \{\beta_j | 1 \leq j\leq N\}$.
Hence $\Delta^x=(\Delta^{\re})^x$. \qed
\end{proposition}

\subsubsection{GRS of Nichols algebras with finite sets of roots}

Let $V$ be a braided vector space of diagonal type with matrix $(q_{jk})$. Let $\chi:\zt\times\zt\to\ku^\times$ be the bicharacter such that $\chi(\alpha_j,\alpha_k)=q_{jk}$ for all $j,k\in\I$.
The set $\Delta^{V}_+$ of degrees of a PBW basis of $\cB(V)$, counted with their multiplicities as in \cite{H-Weyl grp}, does not depend on the PBW
basis. For each $i\in \I$ let $c_{ii}^{V}=2$,
\begin{align}
-c_{ij}^{V}&:= \min \left\{ n \in \N_0: (n+1)_{q_{ii}}
(1-q_{ii}^n q_{ij}q_{ji} )=0 \right\} & & \notag \\
&=\max\{ n\in\N_0: n\alpha_i+\alpha_j\in \Delta_+^{V}\}, & j & \neq i \label{defn:mij}
\end{align}
$s_i^{V}\in\Aut(\Z^\theta)$ such that $s_i^{V}(\alpha_j)=\alpha_j-c_{ij}^{V}\alpha_i$. Let $\cX$ be the set of braided vector spaces $(V,c)$ of diagonal type such that $\Delta_+^{V}$ is finite.
Set $\Delta^{V}=\Delta^{V}_+\cup (-\Delta^{V}_+)$. Define $\rho_i(V)$ as the braided vector space with associated bicharacter $\rho_i(\chi)(\alpha,\beta)= \chi(s_i^{V}(\alpha),s_i^{V}(\beta))$, $\alpha, \beta \in\zt$. $\cR:= \cR(\cX, \rho, \cC, (\Delta^x)_{x\in\cX})$ is the GRS attached to $\cX$. Indeed
$\Delta^{\rho_i(V)}_+ = s_i\left( \Delta^{V}_+ \setminus\{\alpha_i\}\right)
\cup \{\alpha_i\}$ by
\cite{H-Weyl grp}, so \eqref{eq:def root system 3} holds, and \cite[Theorems 6.2, 6.9]{HS-coideal subalg} completes the proof.

\subsection[Lusztig Isomorphisms]{Lusztig Isomorphisms of Nichols algebras of diagonal type}\label{section:Lusztig isomorphisms}
We recall now Lusztig type isomorphisms \cite{H-isom} of Hopf algebras related with Nichols algebras of diagonal type. They can be thought of as generalizations of the isomorphisms of quantized enveloping algebras in \cite{L-libro}. 
\smallbreak

Fix a bicharacter $\chi : \Z^{\theta} \times \Z^{\theta} \to \ku^{\times}$, $q_{ij}=\chi(\alpha_i, \alpha_j)$ and $(V,c)$ the braided vector space of diagonal type with braiding matrix $(q_{ij})$ for a basis $(E_i)_{1\le i\le\theta}$. \emph{From now on we assume that $V\in\cX$; in particular all the $c_{ij}^{V}$ do exist}. Set $(V,c)$ as an element of $\ydzt$ as in \eqref{eq:YD diagonal}, so $T(V)$ is a Hopf algebra in this category. We recall definitions and results from \cite[Section 4.1]{H-isom}.

Let $\cU(V)$ be the algebra generated by elements $E_i$,
$F_i$, $K_i$, $K_i^{-1}$, $L_i$, $L_i^{-1}$, $1 \leq i \leq \theta$,
and relations
\begin{align*}
XY&=YX, & X,Y \in & \{ K_i^{\pm}, L_i^{\pm}: 1 \leq i \leq \theta \}, \\
K_iK_i^{-1}&=L_iL_i^{-1}=1,  & E_iF_j-F_jE_i&=\delta_{i,j}(K_i-L_i).
\\ K_iE_j&=q_{ij}E_jK_i,  & L_iE_j&=q_{ji}^{-1}E_jL_i,
\\ K_iF_j&=q_{ij}^{-1}F_jK_i, & L_iF_j&=q_{ji}F_jL_i.
\end{align*}
If $(W,c)$ denotes the braided vector space of diagonal type corresponding to the transpose of the matrix $(q_{ij})_{i,j\in\I}$,
then $\cU(V)$ is the quantum double of $T(V)\#\ku\Z^{\theta}$ and $T(W)\#\ku\Z^{\theta}$,
see \cite[Proposition 4.6]{H-isom}.
$\cU(V)$ admits a Hopf algebra structure, where the comultiplication satisfies
\begin{align*}
\Delta(K_i^{\pm1})&=K_i^{\pm1} \ot K_i^{\pm1}, & \Delta(E_i)&=E_i \ot 1 + K_i \ot E_i,
\\ \Delta(L_i^{\pm1})&=L_i^{\pm1} \ot L_i^{\pm1}, & \Delta(F_i)&=F_i \ot L_i + 1 \ot F_i.
\end{align*}
Then $\cU(V)$ is a $\Z^{\theta}$-graded Hopf algebra such that
\begin{align*}
\deg(K_i)&=\deg(L_i)=0, & \deg(E_i)&=\alpha_i, & \deg(F_i)&=-\alpha_i.
\end{align*}

We fix the following notation:
\begin{itemize}
\item $\cU^+(V)$ (respectively, $\cU^-(V)$) is the subalgebra
generated by $E_i$ (respectively, $F_i$), $1 \leq i \leq \theta$;
\item $\cU^{\ge0}(V)$ (respectively, $\cU^{\le0}(V)$) is the subalgebra
generated by $E_i$, $K_i$, $K_i^{-1}$ (respectively, $F_i$, $L_i$, $L_i^{-1}$), $1
\leq i \leq \theta$,
\item $\cU^{+0}(V)$ (respectively, $\cU^{-0}(V)$) is the subalgebra
generated by $K_i$, $K_i^{-1}$ (respectively, $L_i$, $L_i^{-1}$), $1
\leq i \leq \theta$, which are isomorphic to $\ku \Z^{\theta}$ as Hopf algebras;
\item $\cU^0(V)$ is the subalgebra
generated by $K_i$, $K_i^{-1}$, $L_i$ and $L_i^{-1}$, which is isomorphic to $\ku \Z^{2\theta}$ as Hopf algebras.
\end{itemize}
Then $\cU^+(V)$ is isomorphic to $T(V)$ and $\cU^{\ge 0}(V)$ is isomorphic to $T(V)\#\ku\Z^{\theta}$ as Hopf algebras.
For each homogeneous element $E\in\cU^+(V)_n$ and $k\in\{0,1,\ldots,n\}$,
$\underline\Delta_{n-k,k}(E)$ is the component of $\underline \Delta(E)$ in $\cU^+(V)_{n-k}\ot\cU^+(V)_k$.

We consider some skew-derivations as in \cite[4.2]{H-isom}.

\begin{proposition}\label{prop:derivadas torcidas}
There exist linear endomorphisms $\partial^K_i$, $\partial^L_i$ of $\cU^+(V)$, $i\in\I$, such that for each $E\in \cU^+(V)_n$, $n\in\N$,
\begin{align*}
\underline \Delta_{n-1,1}(E)&=\sum_{i=1}^\theta \partial_i^K(E) \ot E_i,  & \underline \Delta_{1,n-1}(E)&=\sum_{i=1}^\theta E_i \ot\partial_i^L(E),
\end{align*}
Moreover, $EF_i-F_iE = \partial^K_i(E)K_i- L_i \partial^L_i(E)$ and
\begin{align*}
\partial_i^K(1)&= \partial_i^L(1)=0, & \partial_i^K(EE')&= \partial_i^K(E)(K_i \cdot E')+ E \partial_i^K(E'),
\\ \partial_i^K(E_j)&= \partial_i^L(E_j)=\delta_{i,j}, & \partial_i^L(EE')&= \partial_i^L(E)E'+ (L_i^{-1}\cdot E) \partial_i^L(E').
\end{align*}
for all $E,E' \in \cU^+(V)$, $j\in\I$. \qed
\end{proposition}

\begin{corollary}\label{coro:FE=EF+cosas}
Let $m,n\in\N$, $E\in\cU^+(V)_m$, $F\in\cU^-(V)_n$.
\smallbreak

\noindent\vi If $m\geq n$, then $FE\in EF+\sum_{j=0}^{n-1} \cU^+(V)_{m-n+j}\cU^0(V)\cU^-(V)_j$.
\smallbreak

\noindent\vii If $m<n$, then $FE\in EF+\sum_{j=0}^{m-1} \cU^+(V)_{j}\cU^0(V)\cU^-(V)_{n-m+j}$.
\end{corollary}
\pf
\vi We argue by induction on $n$. If $n=1$, then $F$ is a linear combination of $F_i$'s and Proposition \ref{prop:derivadas torcidas} applies.
Let $n>1$. We may assume that $F=F_iF'$, $F'\in\cU^+(V)_{n-1}$. By inductive hypothesis,
$$ FE=F_iF'E\in F_iEF'+\sum_{j=0}^{n-2} F_i\cU^+(V)_{m-n+1+j}\cU^0(V)\cU^-(V)_j,$$
and the proof follows by Proposition \ref{prop:derivadas torcidas}.
\smallbreak

\noindent\vii We may assume that $F=F_{i_1}\dots F_{i_{n-m}} F'$, $F'\in\cU^+(V)_m$. By \vi,
\begin{align*}
FE=F_{i_1}\dots F_{i_{n-m}}F'E\in & F_{i_1}\dots F_{i_{n-m}}EF'\\
&+\sum_{j=0}^{m-1} F_{i_1}\dots F_{i_{n-m}} \cU^+(V)_{j}\cU^0(V)\cU^-(V)_j,
\end{align*}
and the proof follows again by Proposition \ref{prop:derivadas torcidas}.
\epf

\bigbreak

Fix $i\in\I$. For each $j\neq i$ we set as in \cite{H-isom}:
\begin{align*}
E_{j,0(i)}^+=E_{j,0(i)}^-&=E_j, & F_{j,0(i)}^+=F_{j,0(i)}^-&=F_j,
\end{align*}
and recursively,
\begin{gather}\label{eq:def Ejm+,Fjm+}
\begin{aligned}
E_{j,m+1(i)}^+ &:= E_iE_{j,m(i)}^+ - (K_i \cdot E_{j,m(i)}^+)E_i = (\ad_c E_i)^{m+1}E_j, \\
F_{j,m+1(i)}^+ &:= F_iF_{j,m(i)}^+ - (L_i \cdot F_{j,m(i)}^+)F_i= (\ad_c F_i)^{m+1}F_j,
\\ E_{j,m+1(i)}^- &:= E_iE_{j,m(i)}^- - (L_i \cdot E_{j,m(i)}^-)E_i,
\\ F_{j,m+1(i)}^- &:= F_iF_{j,m(i)}^- - (K_i \cdot F_{j,m(i)}^-)F_i.
\end{aligned}
\end{gather}
When $i$ is explicit, we simply denote $E_{j,m(i)}^\pm$ by $E_{j,m}^\pm$.

\begin{remark} \cite[Lemma 4.23, Corollary 4.25]{H-isom}
For all $N\in\N$,
\begin{align}
E_{j,N}^+ &= \sum_{s=0}^N (-1)^s q_{ij}^s q_{ii}^{s(s-1)/2} \binom{N}{s}_{q_{ii}} E_i^{N-s}E_jE_i^s,
\label{eq:expresion adjunta} \\
E_i^N F_i-F_iE_i^N&= (N)_{q_{ii}}(q_{ii}^{1-N}K_i-L_i)E_i^{N-1}, \label{eq:simple power root contra Fp} \\
E_{j,N}^+ F_i - F_i E_{j,N}^+&= (N)_{q_{ii}}(q_{ii}^{N-1}\widetilde{q_{ij}}-1)L_i E_{j,N-1}^+. \label{eqn:corchete E+ con Fp}
\end{align}
The coproduct satisfies the following identities:
\begin{align}
\underline\Delta \left( E_i^N \right) &= \sum_{s=0}^N \binom{N}{s}_{q_{ii}} E_i^{s} \otimes E_{i}^{N-s},
\label{eq:coproducto Ep a la N} \\
\underline\Delta \left( E_{j,N}^+ \right) &= E_{j,N}^+ \ot 1 + \sum_{s=0}^{N-1} \left( \prod_{r=N-s}^{N-1} (r)_{q_{ii}}(1-q_{ii}^{r}\widetilde{q_{ij}})  \right) E_i^{s} \otimes E_{j,N-s}^+,
\label{eq:coproducto Ei+} \\
\underline\Delta \left( E_{j,N}^- \right) &= 1 \ot E_{j,N}^- \notag \\
& \quad + \sum_{s=0}^{N-1} q_{ij}^s\left( \prod_{r=N-s}^{N-1} (r)_{q_{ii}}\left(1-q_{ii}^{-r}\widetilde{q_{ij}}^{-1}\right)  \right) E_{j,N-s}^- \otimes E_i^{s}.
\label{eq:coproducto Ei-}
\end{align}
\end{remark}

Let $\cJ^{\pm}(V)$ be the ideal of $\cU^\pm(V)$ such that $\u^+(V)=\cU^+(V)/\cJ^+(V)$, $\u^-(V)=\cU^-(V)/\cJ^-(V)$ are isomorphic to $\cB(V)$, $\cB(W)$, respectively, and set $\u(V)= \cU(V) / (\cJ^-(V)+\cJ^+(V))$.
Then $\u(V)$ is the quantum double of the algebras $\cB(V)\#\ku\Z^\theta$ and $\cB(W)\#\ku\Z^\theta$ by \cite[Theorem 5.8]{H-isom}.

We need another quotients of $\cU(V)$ in order to introduce Lusztig isomorphisms. First we recall
\cite[Definition 2.6]{A-presentation}.
An element $i\in\I$ is a \emph{Cartan vertex} of $V$ if $\widetilde{q_{ij}}=q_{ii}^{c_{ij}^{V}}$ for all $j \neq i$.
The set of \emph{Cartan roots} is
\begin{equation}\label{eq:def O(chi)}
\cO(V):=\{ s_{i_1}^{V}\dots s_{i_k}(\alpha_i): i \mbox{ is a Cartan vertex of } \rho_{i_k}\dots\rho_{i_1}(V)\}.
\end{equation}
Set $N_i= \ord q_{ii}$.
$\cJ_i^+(V)$, $\cJ_i^-(V)$ are the ideals of $\cU^+(V)$, respectively $\cU^-(V)$, generated by
\begin{itemize}
  \item[(a)] $E_i^{N_i}$, respectively $F_i^{N_i}$, if $i$ is not a Cartan vertex,
  \item[(b)] $E_{j,-c_{ij}^{V}+1}^+$, respectively $F_{j,-c_{ij}^{V}+1}^+$, for each $i$ such that $N_i\geq-c_{ij}^{V}+1$.
\end{itemize}
Set also
\begin{align*}
\cU_i(V)&:=\cU(V)/\left(\cJ_i^+(V)+\cJ_i^-(V) \right), & \cU_i^\pm(V)&:=\cU^\pm(V)/\cJ_i^\pm(V).
\end{align*}

Let $\Eb_j$, $\Fb_j$, $\Kb_j$, $\Lb_j$ be the generators of $\cU(\schi)$. Set
\begin{equation}\label{eqn:escalares lambda}
\lambda_j(V):= (-c_{ij}^{V})_{q_{ii}}\prod_{s=0}^{-c_{ij}^{V}-1}(q_{ii}^s\widetilde{q_{ij}}-1) ,  \qquad j \neq i.
\end{equation}

\begin{theorem}{\cite[Lemma 6.5, Theorem 6.12]{H-isom}}\label{thm: iso Lusztig-Heck}
There exist algebra maps
\begin{equation}\label{eqn:iso lusztig para cUp}
T_i, T_i^-: \cU(V) \to \cU_i(\schi)
\end{equation}
univocally determined by the following conditions:
\begin{align*}
T_i(K_i)&=T_i^-(K_i)=\Kb_i^{-1}, & T_i(K_j)&=T_i^-(K_j)=\Kb_i^{-c_{ij}^{V}}\Kb_j,
\\ T_i(L_i)&=T_i^-(L_i)=\Lb_i^{-1}, & T_i(L_j)&=T_i^-(L_j)=\Lb_i^{-c_{ij}^{V}}\Lb_j,
\\ T_i(E_i)&=\Fb_i\Lb_i^{-1}, & T_i(E_j)&=\Eb^+_{j,-c_{ij}^{V}},
\\ T_i(F_i)&=\Kb_i^{-1}\Eb_i, & T_i(F_j)&=\lambda_i(\schi)^{-1}\Fb^+_{j,-c_{ij}^{V}},
\\ T_i^-(E_i)&=\Kb_i^{-1}\Fb_i, & T_i^-(E_j)&=\lambda_i(\schi^{-1})^{-1}\Eb^-_{j,-c_{ij}^{V}},
\\ T_i^-(F_i)&=\Eb_i\Lb_i^{-1}, & T_i^-(F_j)&=\Fb^-_{j,-c_{ij}^{V}}.
\end{align*}
Such morphisms induce algebra isomorphisms (denoted by the same name):
\begin{align}
T_i, T_i^-: \, &\u(V) \to \u(\rho_i(V)) &  \mbox{such that }& T_iT_i^-=T_i^-T_i=\id.\qed
\end{align}
\end{theorem}


\subsection{Lusztig isomorphisms and PBW bases}\label{subsection:pbw}

Fix an reduced expression $w=\sigma_{i_1}^{V} \sigma_{i_2}\cdots \sigma_{i_M}$ of the element of maximal length of
$\cW(V)$. If $1\leq k\leq M$ then set $\beta_k:= s_{i_1}^{V}\cdots s_{i_{k-1}}(\alpha_{i_k})$,
$q_k:=\chi(\beta_k,\beta_k)$, $N_k = \ord q_k\in\N\cup\{\infty\}$.
By Proposition \ref{Prop:maxlongitudCH} $\beta_k\neq\beta_l$ if $k\neq l$, and $\Delta_+^{V}=\{\beta_k|1\leq k\leq M\}$.
Let
\begin{align}\label{eq:PBW generators}
E_{\beta_k}&=T_{i_1}\cdots T_{i_{k-1}}(E_{i_k})\in\u(V)^+_{\beta_k}, & F_{\beta_k}&=T_{i_1}\cdots T_{i_{k-1}}(F_{i_k})\in\u(V)^-_{\beta_k}.
\end{align}
and for each $\ab=(a_1,\dots,a_M)\in\N_0^M$,
\begin{align}\label{eq:def E a la a, F a la a}
\E^{\ab}&=E_{\beta_M}^{a_M}E_{\beta_{M-1}}^{a_{M-1}} \cdots
E_{\beta_1}^{a_1}, & \F^{\ab}&= F_{\beta_M}^{a_M}F_{\beta_{M-1}}^{a_{M-1}} \cdots
F_{\beta_1}^{a_1}.
\end{align}

\begin{theorem} \label{thm: HY PBW bases}
\cite[Theorems 4.5, 4.8, 4.9]{HY-shapov}
The sets
\begin{align*}
& \{ \E^{\ab}\, | \,  \ab\in\N_0^M, \, 0\leq a_k < N_k, \, 1\leq k\leq
M\},
\\ & \{ \F^{\ab} \, | \, \ab\in\N_0^M, \, 0\leq a_k < N_k, \, 1\leq k\leq
M\},
\end{align*}
are bases of the vector spaces $\u^+(V)$, $\u^-(V)$, respectively.\qed
\end{theorem}


\section[Distinguished pre-Nichols algebras]{Distinguished pre-Nichols algebras}
\label{section:Distinguished pre-Nichols algebras}

From now on we assume that $\cB(V)$ \emph{is finite-dimensional}.

\subsection{\hspace{1pt}} We now consider some intermediate quotients between $\cU(V)$ and $\u(V)$.

\begin{definition}
Let $\cI(V)$ be the ideal of $T(V)$ generated by all
the relations in \cite[Theorem 3.1]{A-presentation}, except the power root vectors $E_\alpha^{N_\alpha}$, $\alpha\in\cO(V)$,
plus the quantum Serre relations $(\ad_c E_i)^{1-c_{ij}^{V}} E_j$ for those $i\neq j$ such that
$q_{ii}^{c_{ij}^{V}}=q_{ij}q_{ji}=q_{ii}$.
Then $\dpn(V)=T(V)/\cI(V)$ is the \emph{distinguished pre-Nichols algebra} of $(V,c)$.
\end{definition}

We identify $\cI(V)$ as an ideal $\cI^+(V)$ of $\cU^+(V)$. Let $\cI^-(V)$ be the corresponding ideal of $\cU^-(V)$.
We denote by $U(V)$ the quotient of $\cU(V)$ by the ideal generated by $\cI^+(V)$ and $\cI^-(V)$. By abuse of notation we denote by $E_i$, $F_i$, $K_i^\pm$, $L_i^{\pm}$ the generators of $U(V)$.
$U^+(V)=\dpn(V)$, $U^-(V)$ are, respectively, the subalgebras of $U(V)$ generated by $E_i$, $F_i$. Let $U^0(V)$ be the subalgebra generated by $K_i^\pm$, $L_i^{\pm}$.

\begin{proposition}\cite[Proposition 3.3]{A-presentation}\label{prop:proyeccion sobre Nichols}
$U(V)$ is a Hopf algebra and there exist
a canonical Hopf algebra morphism  $\pi_V:U(V)\twoheadrightarrow\u(V)$ such that
$\pi_V\left(U^\pm(V)\right)=\u^\pm(V)$. The multiplication $ m:U^+(V)\ot U^0(V)\ot U^-(V)\to U(V) $
gives an isomorphism of graded vector spaces. \qed
\end{proposition}

$U(V)$ is the quantum double of $\dpn(V)\# \ku\zt$ and $\dpn(W)\# \ku\zt$ since there exists a Hopf pairing induced by the one between $T(V)\# \ku\zt$ and $T(W)\# \ku\zt$.

\begin{remark}
Let $\g$ be a finite-dimensional semisimple Lie algebra over $\C$, $C=(c_{ij})\in\Z^{\I\times\I}$ its finite Cartan matrix and $D=\operatorname{diag}(d_1,\dots,d_\theta)$ such that $DA$ is symmetric. If $q$ is a root of unity of odd order, then the symmetric matrix $(q_{ij})_{i,j\in\I}$, $q_{ij}=q^{d_ic_{ij}}$, defines a braiding of Cartan type. The small quantum group $\u_q(\g)$ is the quotient of $\u(V)$ by the central elements $K_i-L_i^{-1}$, $1\le i\le\theta$, while the quantized enveloping algebra $U_q(\g)$ is obtained from $U(V)$ analogously.
\end{remark}

The Lusztig isomorphisms descend to the family of algebras $U(V)$.

\begin{proposition}\cite[Proposition 3.26]{A-presentation}\label{prop:iso Lusztig para las U}
The maps \eqref{eqn:iso lusztig para cUp} induce algebra isomorphisms $ T_i, T_i^-: U(V) \to U(\schi)$
such that $T_i T_i^-=T_i^- T_i= \id_{U(V)}$.\qed
\end{proposition}

\begin{remark}\label{rem:PRV=0 si alpha no cartan}
Using the Lusztig isomorphisms we deduce that
\begin{align*}
E_\alpha^{N_\alpha}, F_\alpha^{N_\alpha}&\neq0 \mbox{ for all }\alpha\in\cO(V), & E_\alpha^{N_\alpha}&=F_\alpha^{N_\alpha}=0 \mbox{ for all }\alpha\notin\cO(V),
\end{align*}
on $U(V)$. The Hilbert series of $U^\pm(V)$ is:
$$ \cH_{U^\pm(V)}= \left( \prod_{\alpha\in \Delta_+^{V}\setminus\cO(V) } (t^{\alpha})_{N_\alpha}\right)\bigcup \left( \prod_{\alpha\in \cO(V) } \frac{1}{1-t^{\alpha}} \right).  $$
Indeed $U^+(V)$ has a PBW basis of Lyndon hyperwords as in \cite{Kh} with the same PBW generators of $\u^+(V)$, see the proof of \cite[Theorem 3.1]{A-presentation}.
\end{remark}

Set $E_\alpha, \E^{\ab}\in U^+(V)$, $F_\alpha, \F^{\ab}\in U^-(V)$, $\alpha\in\Delta_+^{V}$, $\ab\in\N_0^M$, as in \eqref{eq:PBW generators}, \eqref{eq:def E a la a, F a la a}.

\begin{theorem} \label{thm: PBW bases preNichols}
The sets
\begin{align*}
& \{ \E^{\ab}\, | \,  \ab\in\N_0^M, \, 0\leq a_k < N_k \mbox{ if }\beta_k\notin\cO(V) \}, \\
& \{ \F^{\ab}\, | \,  \ab\in\N_0^M, \, 0\leq a_k < N_k \mbox{ if }\beta_k\notin\cO(V) \},
\end{align*}
determine bases of the vector space $U^+(V)$, $U^-(V)$, respectively.
\end{theorem}
\pf
The expression of the Hilbert series reduces the problem to the linearly independence of this set, which is proved following the same recursion as \cite[Theorem 4.5]{HY-shapov}.
\epf

\subsection{\hspace{1pt}}
Let $\cB_{i}$ be the algebra generated by $E_i$ on $U(V)$. This is a braided graded Hopf algebra. Its graded dual $\cB_{i}^*$ is also a graded braided Hopf algebra.
There exists a projection $\pi_{i,V}:U^+(V)\twoheadrightarrow\cB_{i}$ of braided graded Hopf algebras annihilating all the $E_j$ for $j\neq i$; the inclusion $\iota_{i,V}:\cB_{i}\hookrightarrow U^+(V)$ is a section for this projection.

\begin{remark}\label{rem:dual algebra vertice cartan}
Let $i$ be a Cartan vertex. $\cB_{i}^*$ has a basis $\{\Ef_i^{(n)}:n\in\N_0\}$, where $\Ef_i^{(n)}(E_i^m)=\delta_{n,m}$. The algebra and the coalgebra structures satisfy:
\begin{align}\label{eq:producto coproducto dual}
\Ef_i^{(j)}\cdot \Ef_i^{(k)}&= \binom{j+k}{j}_{q_{ii}}\Ef_i^{(k+j
)}, & \underline{\Delta}(\Ef_i^{(n)})&=\sum_{j=0}^n \Ef_i^{(j)}\otimes \Ef_i^{(n-j)},
\end{align}
with unit $\Ef_i^{(0)}=1$, and counit $\eps(\Ef_i^{(n)})=\delta_{n,0}$. In particular $\Ef_i^{(j)}\cdot \Ef_i^{(N_i-j)}=0$ if $1\leq j\leq N_i-1$. As an algebra $\cB_{i}^*$ is generated by $\Ef_i^{(1)}$, $\Ef_i^{(N_i)}$.
\end{remark}

\begin{remark}\label{rem:accion dual Bi}
There exist left and right actions of $\cB_{i}^*$ on $U^+(V)$ given by
\begin{align}\label{eq:acciones dual}
\Ef\triangleright Y & = Y_{(1)} \Ef (\pi_{i,V}(Y_{(2)})), & Y \triangleleft\Ef & =  \Ef(\pi_{i,V}(Y_{(1)}))Y_{(2)},
\end{align}
for each $\Ef\in \cB_{i}^*$, $Y\in U^+(V)$.
In particular,
\begin{align}\label{eq:accion X1}
\Ef_i^{(1)}\tri Y & = \partial_i^K(Y), &  Y \triangleleft \Ef_i^{(1)} & =\partial_i^L(Y), & \mbox{for all }& Y\in U^+(V).
\end{align}
\end{remark}

\smallbreak

\begin{lemma}
For all $X\in U^+(V)_\beta$, $Y\in U^+(V)$, $t\in\N$,
\begin{equation}\label{eq:Efi como skew derivation}
(XY)\triangleleft\Ef_i^{(t)}= \sum_{r=0}^t \chi(\beta-r\alpha_i,\alpha_i)^{t-r} \left(X\triangleleft\Ef_i^{(r)}\right) \left(Y\triangleleft\Ef_i^{(t-r)}\right).
\end{equation}
\end{lemma}
\pf
Notice that
\begin{align*}
\sum_{t\geq 0} E_i^t  \ot (XY)\triangleleft\Ef_i^{(t)} & =(\pi_{i,V}\ot\id)\underline\Delta(XY)= (\pi_{i,V}\ot\id)\left(\underline\Delta(X)
\cdot \underline\Delta(Y)\right) \\
& = \left( \sum_{r\geq 0} E_i^r\ot X\triangleleft\Ef_i^{(r)}  \right) \cdot \left(  \sum_{s\geq 0} E_i^s\ot Y\triangleleft\Ef_i^{(s)} \right) \\
& = \sum_{r,s\geq 0} \chi(\beta-r\alpha_i,\alpha_i)^s \, E_i^{r+s} \ot (X\triangleleft\Ef_i^{(r)}) ( Y\triangleleft\Ef_i^{(s)}).
\end{align*}
Then compare the terms of the form $E_i^t\ot -$.
\epf

\subsection{\hspace{1pt}}
Let $U^+_{\pm i}(V)$ be the subalgebra generated by $E_{j,N}^\pm$, $j\neq i$, $N\in\N_0$.
Given $\alpha=\sum\limits_{i=1}^\theta n_i\alpha_i\in\zt$, set $K_\alpha=\prod\limits_{i=1}^\theta K_i^{n_i}$, $L_\alpha=\prod\limits_{i=1}^\theta L_i^{n_i}\in U^0(V)$.

\begin{lemma}\label{lema:kernel derivaciones generado por ad}
\vi $U_{+ i}^+(V)=U^+(V)^{\co \pi_{i,V}}$, $U_{-i}^+(V)=^{\co \pi_{i,V}}U^+(V)$. Then there exist isomorphisms of graded vector spaces $U^+(V) \cong U_{\pm i}^+(V) \ot \cB_{i}$.

\smallbreak
\noindent\vii If $i$ is not a Cartan vertex, then $\ker(\partial_i^K)=U_{+ i}^+(V)$, $\ker(\partial_i^L)=U_{-i}^+(V)$.

\smallbreak
\noindent\viii If $i$ is a Cartan vertex, then $\ker(\partial_i^K)= U_{+ i}^+(V)\ku\left[E_i^{N_i}\right]$,
$\ker(\partial_i^L)= U_{-i}^+(V)\ku\left[E_i^{N_i}\right]$.
\end{lemma}
\pf
The claims about $U_{+i}^+(V)$ follow by \cite[Lemma 2.4]{A-standard}, \cite[Lemma 4.31]{H-isom}, so we prove those
about $U_{-i}^+(V)$.
By \eqref{eq:coproducto Ei-} $U_{-i}^+(V)$ is a right coideal subalgebra contained in $B=^{\co \pi_{i,V}}U^+(V)$. We claim that $U_{-i}^+(V) \cB_{i}$ is a left ideal of $U^+(V)$. Indeed if $E\in U_{-i}^+(V)$, then $E_jE\in U_{-i}^+(V)$ for $j\neq i$ since $E_j\in U_{-i}^+(V)$, and $E_i X\in U_{-i}^+(V) \cB_{i}$ since
\begin{align*}
E_i E_{j,n}^- &= E_{j,n+1}^- + \chi(n\alpha_i+\alpha_j,-\alpha_i) E_{j,n}^- E_i, & j\neq i, \, & n\in\N_0.
\end{align*}
But $1\in U_{-i}^+(V) \cB_{i}$, so $U_{-i}^+(V) \cB_{i}=U^+(V)$. By Proposition \ref{pro:masuoka} \viii, there exists an isomorphism $B\otimes \cB_{i}\simeq U^+(V)$, so $B= U_{-i}^+(V)$ and the multiplication gives this isomorphism $U_{-i}^+(V) \ot \cB_{i} \simeq U^+(V)$. Therefore for each $X\in U^+(V)$ there exist unique $X_n\in U_{-i}^+(V)$ such that $X= \sum_{n\geq 0}X_n E_i^n$. If $X$ has degree $\beta$, then $\partial_i^L(X) = \sum_{n\geq 1} (n)_{q_{ii}} \chi(n\alpha_i-\beta,\alpha_i) \, X_n  E_i^{n-1}$ since $U_{-i}^+(V)\subseteq \ker\partial_i^L$. If $X\in\ker\partial_i^L$, then $X_n=0$ for all $n\notin\N N_i$.
\epf

\medbreak

\begin{remark}\label{rem:algebra de Hopf YD para U+i}
$U_{+ i}^+(V)$ is a braided Hopf algebra in $^{\cB_{i}\#\ku\Z^\theta}_{\cB_{i}\#\ku\Z^\theta}\mathcal{YD}$. Indeed Lemma \ref{lema:kernel derivaciones generado por ad}\viii says that $U_{+ i}^+(V)=U^+(V)^{\co \pi_{i,V}}$ , and by \cite[Lemma 3.1]{AHS}
$$ U^+(V)^{\co \pi_{i,V}}=(U^+(V)\#\ku\Z^\theta)^{\co \pi_{i,V}\#\id} .$$
The action and coaction satisfy, for each $E\in U^+_{+i}(V)_\beta$,
\begin{align}\label{eq:accion-coaccion a izquierda}
E_i\rightharpoonup E&=(\ad_c E_i) E, & \lambda(E)&= \sum_{n\geq 0} E_i^n K_{\beta-n\alpha_i}\ot E\triangleleft \Ef_i^{(n)}.
\end{align}
\end{remark}

\begin{remark}\label{rem:algebra de Hopf YD para U-i}
If $R$ is braided Hopf algebra $R$ in $\ydzt$, then there exists an structure of braided Hopf algebra $R^{\bop}$ as in \cite[Proposition 2.2.4]{AG} with underlying Yetter-Drinfeld module $R$ and
\begin{align*}
m^{\bop}&=m\circ c_{R,R} , & \underline{\Delta}^{\bop}&= c_{R,R}^{-1}\circ \underline{\Delta}, & \mathcal{S}^{\bop}&= \mathcal{S}.
\end{align*}
This applies for $U^+(V)$, $\cB_{i}$, and $\pi_{i,V}:U^+(V)^{\bop}\twoheadrightarrow\cB_{i}^{\bop}$ is a Hopf algebra map in $\ydzt$. Consider the Hopf algebras $U^+(V)^{\bop}\#\ku\Z^\theta$, $\cB_{i}^{\bop}\#\ku\Z^\theta$, and the Hopf algebra maps $\pi_{i,V}\#\id$, $\iota_{i,V}\#\id$. Then
$$ (U^+(V)^{\bop}\#\ku\Z^\theta)^{\co \pi_{i,V}\#\id}= (U^+(V)^{\bop})^{\co \pi_{i,V}}=^{\co \pi_{i,V}}U^+(V)=U_{-i}^+(V) $$
is a braided Hopf algebra in $^{\cB_{i}^{\bop}\#\ku\Z^\theta}_{\cB_{i}^{\bop}\#\ku\Z^\theta}\mathcal{YD}$.
The action and coaction satisfy
\begin{align}\label{eq:accion-coaccion para U-i}
E_i\rightharpoonup' E&=[E,E_i]_c, & \lambda'(E)&= \sum_{n\geq 0} E_i^n K_{\beta-n\alpha_i}\ot \Ef_i^{(n)} \triangleright E,
\end{align}
for each $E\in U^+_{-i}(V)_\beta$. But $U_{-i}^+(V)$ has the opposite product, so we take the braided Hopf algebra structure on $U_{-i}^+(V)$ obtained by applying the $\bop$ construction. The coproduct $\Delta_i:U^+_{-i}(V)\to U^+_{-i}(V)\underline{\ot} U^+_{-i}(V)$ is given by
\begin{align}\label{eq:formula delta sub i}
\Delta_i(x) &= x_{(1)} \ot \iota_{i,V}\pi_{i,V}( \mathcal{S}^{-1}(x_{(2)})) x_{(3)} , &   x\in & U^+_{-i}(V)
\end{align}
Here $U^+_{-i}(V)\underline{\ot} U^+_{-i}(V)$ denotes the space $U^+_{-i}(V)\underline{\ot} U^+_{-i}(V)$ with the algebra structure viewed as an element in $^{\cB_{i}^{\bop}\#\ku\Z^\theta}_{\cB_{i}^{\bop}\#\ku\Z^\theta}\mathcal{YD}$,
\begin{align*}
(x\ot y) \cdot (w\ot z)
& = \sum_{n=0}^{N_i-1} \chi(\beta,\gamma-n\alpha_i)
\,  x \left(\Ef_i^{(n)} \triangleright w \right) \underline{\ot} \left(E_i^n \rightharpoonup' y  \right) z,
\end{align*}
for $x,y,w,z\in U_{-i}^+(V)$, where $y$, $w$ are homogeneous of degrees $\beta,\gamma\in\N_0^\theta$.
\end{remark}

\begin{remark}\label{rem:Ti(U-)=U+} By \cite[Lemma 6.7]{H-isom},
\begin{align}
q_{ii}\partial_i^L \left(T_i(X)\right) & 
=-\chi(\beta,\alpha_i)^{-1} T_i(E_i \rightharpoonup' X),
\label{eq:Ti con delta L}\\
T_i \left(K_i^{-1}\cdot\partial_i^K(X)\right) & 
= -\left( \Eb_\rightharpoonup T_i(X)\right),
\label{eq:Ti con delta K}\\
T_i(E_{j,n}^-)&= \underline{q}_{ii}^{\,\, n} \left(\prod_{t=-c_{ij}^{V}-n}^{-c_{ij}^{V}-1}
(t+1)_{\underline{q}_{ii}}(1-\underline{q}_{ii}^{\,\, t}\underline{q}_{ij}\underline{q}_{ji}) \right) \Eb_{j,-c_{ij}^{V}-n}^+, \label{eq:Ti(E-)=E+}
\end{align}
for each $X\in U^+_{-i}(V)_\beta$, $\beta\in\N_0^\theta$, $n\in\N_0$, where $(\underline{q}_{jk})_{j,k\in\I}$ is the matrix of $\rho_i(V)$. In particular $T_i(U^+_{-i}(V))= U^+_{+i}(\rho_i(V))$.
\end{remark}

\bigbreak

\subsection{\hspace{1pt}}
Next results resemble \cite[Lemma 4.6, Theorem 4.8]{HY-shapov}.

\begin{lemma}\label{lem:base de ker pi}
The set
\begin{align}\label{eq:base de U+i}
& \{ \E^{\ab}\, | \,  \ab\in\N_0^M, \, a_1=0, \, 0\leq a_k < N_k \mbox{ if }\beta_k\notin\cO(V) \}
\end{align}
is a basis of the vector space $U^+_{+i}(V)$.
\end{lemma}
\pf
Set $i=i_1$. By Lemma \ref{lema:kernel derivaciones generado por ad} \vi, $E_{\beta_\ell}=\sum_{k\geq 0} X_k E_i^k$ for unique $X_k\in\ U^+_{+i}(V)$. By \eqref{eq:Ti(E-)=E+}, $T_i^-(U^+_{+i}(V))=U^+_{-i}(\rho_i(V))$, so $T_i^-(X_k)\in U^+(\rho_i(V))$. As
$$ \sum_{k\geq 0} T_i^-(X_k) (K_i^{-1}F_i)^k=T_i^-(E_{\beta_\ell})=T_{i_2}\dots T_{i_{\ell-1}}(E_{i_l})\in U^+(V),$$
last statement of Proposition \ref{prop:proyeccion sobre Nichols} implies that $X_k=0$ for all $k>0$, so $E_{\beta_\ell}=X_0\in U^+_{+i}(V)$. Then the set \eqref{eq:base de U+i} is contained in $U^+_{+i}(V)$. By Lemma \ref{lema:kernel derivaciones generado por ad} \vi and Theorem \ref{thm: PBW bases preNichols}, this set should also generate $U^+_{+i}(V)$ because it generates a subspace with the same Hilbert series.
\epf

\begin{proposition} \label{prop: corchete entre Ebetas}
For each pair $1\leq k< \ell \leq M$,
\begin{align}
E_{\beta_k}E_{\beta_\ell}- \chi(\beta_k,\beta_\ell)E_{\beta_\ell}E_{\beta_k}
&=\sum c_{a_{k+1},\ldots,a_{\ell-1}}  E_{\beta_{\ell-1}}^{a_{\ell-1}} \cdots E_{\beta_{k+1}}^{a_{k+1}} \in \u^+(V), \label{eq:corchete entre Ebetas}
\\ F_{\beta_k}F_{\beta_\ell}- \chi(\beta_k,\beta_\ell) F_{\beta_\ell}F_{\beta_k}&= \sum d_{a_{k+1},\ldots,a_{\ell-1}}
F_{\beta_{\ell-1}}^{a_{\ell-1}}\cdots F_{\beta_{k+1}}^{a_{k+1}} \in \u^-(V),\label{eq:corchete entre Fbetas}
\end{align}
for some scalars $c_{a_{k+1},\ldots,a_{\ell-1}}, d_{a_{k+1},\ldots,a_{\ell-1}}\in\ku$.
\end{proposition}
\pf
Assume first $k=1$. By Lemma \ref{lem:base de ker pi} there exist $\mathtt{c}_{\ab}\in\ku$ such that
$$ E_{i_1}E_{\beta_\ell}-\chi(\alpha_{i_1},\beta_\ell)E_{\beta_\ell}E_{i_1}= \sum_{\ab: \, a_1=0} \mathtt{c}_{\ab} \, \E^{\ab}. $$
Let $V'=\rho_{i_\ell}\dots \rho_{i_1}V$. By applying $T_{i_\ell}^-\dots T_{i_1}^-$ to the left hand side we obtain an element of $U^{\le 0}(V')$ . Then $\mathtt{c}_{\ab}=0$ if some $a_j>0$ for $j>\ell$ since
\begin{align*}
T_{i_\ell}^-\dots T_{i_1}^-(E_{\beta_j})&=T_{i_\ell}^-\dots T_{i_{j+1}}^-(K_{i_j}^{-1}F_{j})\in U^{\le 0}(V') & \mbox{if } & j\leq \ell, \\
T_{i_\ell}^-\dots T_{i_1}^-(E_{\beta_j})&=T_{i_{\ell+1}}\cdots T_{i_{j-1}}(E_{i_j})\in U^+(V') & \mbox{if } & j> \ell.
\end{align*}

Finally if $\mathtt{c}_{0,a_2,\ldots,a_\ell,0,\ldots, 0}\neq 0$, then
$\sum_{j=2}^\ell a_j\beta_j=\alpha_{i_1}+\beta_\ell$ since $U^+(V)$ is $\zt$-graded so $a_\ell=0$. Therefore we
get \eqref{eq:corchete entre Ebetas} for  $c_{a_{k+1},\ldots,a_{\ell-1}}=\mathtt{c}_{0,a_2,\ldots,a_{\ell-1},0,\ldots, 0}$.
The proof of \eqref{eq:corchete entre Fbetas} is analogous.
\epf

Now we introduce algebra filtrations of $U^+(V)$, $U^{\ge0}(V)$ and $U(V)$ related with the PBW basis.
We order $\N_0^M$, $\N_0^{2M+1}$ lexicographically. In particular $\delta_1<\dots<\delta_M$, where $\{\delta_j\}_{1\le j\le M}$ denotes  the canonical basis of $\Z^M$ to avoid confusion with the basis $\{\alpha_i\}_{1\le i\le\theta}$ of $\zt$.

\smallbreak\noindent\vi
In $U^+(V)$ set $U^+(V)(\ab)$ as the subspace spanned by $\E^{\bb}$, $\bb\le\ab$. This is vector space filtration of $U^+(V)$ so consider the $\N_0^M$-graded vector space
\begin{align*}
\gr U^+(V)&=\oplus_{\ab\in\N_0^M} U^+(V)_{\ab}, & U^+(V)_{\ab}&=U^+(V)(\ab)/\sum_{\bb<\ab} U^+(V)(\bb).
\end{align*}

\smallbreak\noindent\vii
For $U^{\ge0}(V)$, $U^{\ge0}(V)(\ab)$ is the subspace spanned by $\E^{\bb}K^\alpha$, $\bb\le\ab$, $\alpha\in\zt$. In particular $U^{\ge0}(V)(0)=U^{+0}(V)$.

\smallbreak\noindent\viii
For each $\E^{\ab}K^\alpha L^\beta\F^{\bb}\in U(V)$ set
$$ d(\E^{\ab}K^\alpha L^\beta\F^{\bb})=\left(\sum_{j=1}^M (a_j+b_j)\hgt(\beta_j),a_1,\dots,a_M,b_1,\dots,b_M \right) \in\N_0^{2M+1}. $$
For each $\ub\in\N_0^{2M+1}$ let $U(V)(\ub)$ be the subspace spanned by $\E^{\ab}K^\alpha L^\beta\F^{\bb}$, $d(\E^{\ab}K^\alpha L^\beta\F^{\bb})\le\ub$. Then take the associated $\N_0^{2M+1}$-graded vector space
\begin{align*}
\gr U(V)&=\oplus_{\ub\in\N_0^{2M+1}} U(V)_{\ub}, & U(V)_{\ub}&=U(V)(\ub)/\sum_{\vb<\ub} U(V)(\vb).
\end{align*}
Next result generalizes \cite[Proposition 1.7]{DK-root of 1}, \cite[Proposition 10.1]{DP-quantum}

\begin{proposition}\label{prop:filtracion via PBW}
The $\N_0^M$-filtrations on $U^{+}(V)$, $U^{\ge0}(V)$ and the $\N_0^{2M+1}$-filtration on $U(V)$ are algebra filtrations.
\end{proposition}
\pf
First we consider $U^{+}(V)$. We claim that $\E^{\ab}\E^{\bb}\in U^+(V)(\ab+\bb)$. Let $\ab=\delta_k$, $\bb=\delta_\ell$, $1\le k,\ell \le M$.
Note that $U^+(V)(\delta_k)=\ku E_{\beta_k}$, and $E_{\beta_k}E_{\beta_\ell}\in U^+(V)(\delta_k+\delta_\ell)$ if $k\geq \ell$ by definition. Also $E_{\beta_k}E_{\beta_\ell}\in U^+(V)(\delta_k+\delta_\ell)$ if $k<\ell$ by Proposition \ref{prop: corchete entre Ebetas}, since
\begin{align*}
\sum_{j=k+1}^{\ell+1} a_j\delta_j&<\delta_k+\delta_\ell & \mbox{ so }& E_{\beta_{\ell-1}}^{a_{\ell-1}} \cdots E_{\beta_{k+1}}^{a_{k+1}}\in \sum_{\vb<\delta_k+\delta_\ell} U^+(V)(\vb).
\end{align*}
Using this case we can reorder the PBW generators of $\E^{\ab}\E^{\bb}$ for any $\ab,\bb$.

The proof for $U^{\ge0}(V)$ follows from the previous case since $K_i$ $q$-commutes with all the elements of $U^+(V)$. For $U(V)$ we argue as above and reduce the proof to the product between $E_{\beta_k}$, $F_{\beta_\ell}$. By Corollary \ref{coro:FE=EF+cosas},
\begin{align}\label{eq:Fbeta Egamma}
F_{\beta_\ell}E_{\beta_k}\in E_{\beta_k}F_{\beta_\ell}+ \sum_{\vb: v_1\leq \hgt(\beta_\ell)+\hgt(\beta_k)-2}U(V)(\vb)
\end{align}
so
$ F_{\beta_\ell}E_{\beta_k}\in U(V)(\hgt(\beta_\ell)+\hgt(\beta_k),\delta_k,\delta_\ell)$.
\epf

We give a presentation of the corresponding graded algebras.

\begin{corollary}\label{coro:graduada via PBW, presentacion}
\vi The algebra $\gr U^+(V)$ is presented by generators $\Et_k$, $1\le k\le M$ and relations
\begin{align}\label{eq:rels grU+}
\Et_k \Et_\ell&=\chi(\beta_k,\beta_\ell)\, \Et_\ell \Et_k, \quad k<\ell, & \Et_k^{N_k}&=0,\quad \beta_k\notin\cO(V).
\end{align}

\smallbreak
\noindent\vii The algebra $\gr U^{\ge 0}(V)$ is presented by generators $\Et_k$, $1\le k\le M$, $K_i^{\pm}$, $i\in\I$, and relations \eqref{eq:rels grU+},
\begin{align}\label{eq:rels grU ge0 - 1}
K_i K_j &= K_j K_i, & K_iK_i^{-1}&=K_i^{-1}K_i=1, \\
K_i\Et_k &=\chi(\alpha_i,\beta_k)\, \Et_k K_i, &  1\le k\le & M, \, i,j\in\I. \label{eq:rels grU ge0 - 2}
\end{align}
\end{corollary}
\pf
For \vi Let $\mathcal{F}$ be the free algebra generated by $\Et_k$, $1\le k\le M$ and $\pi:\mathcal{F}\to \gr U^+(V)$ the algebra map such that $\pi(\Et_k)=E_{\beta_k}$. By \eqref{eq:corchete entre Ebetas}, $E_{\beta_k}E_{\beta_\ell}= \chi(\beta_k,\beta_\ell)\, E_{\beta_\ell}E_{\beta_k}$ holds in $\gr U^+(V)$, and also $E_{\beta_k}^{N_k}=0$ for $\beta_k\notin\cO(V)$ by Remark \ref{rem:PRV=0 si alpha no cartan}. Then $\pi$ factors through the algebra defined by the relations \eqref{eq:rels grU+}, which is the quotient of a $q$-polynomial algebra of $M$ generators by $\Et_k^{N_k}$ for those generators corresponding to $\beta_k\notin\cO(V)$, so it has a basis
$$ \{\Et_M^{a_m}\dots \Et_1^{a_1} | a_k\in\N_0, a_k<N_k \mbox{ if }\beta_k\notin\cO(V)\}, $$
As $\E^{\ab}\in U^+(V)(\ab)-\sum_{\bb<\ab} U^+(V)(\bb)$, $\gr U^+(V)$ has a basis corresponding to the images of \eqref{eq:base de U+i}. But
$\pi(\Et_M^{a_m}\dots \Et_1^{a_1})=\E^{\ab}$ so $\pi$ is an isomorphism. The proof of \vii follows similarly.
\epf

\begin{corollary}\label{coro:graduada via PBW, presentacion Uchi}
The algebra $\gr U(V)$ is presented by generators $\Et_k$, $\Ft_k$, $1\le k\le M$, $K_i$, $K_i^{-1}$, $L_i$, $L_i^{-1}$, $i\in\I$, and relations
\begin{align*}
XY&=YX, & X,Y \in & \{ K_i^{\pm1}, L_i^{\pm1}: 1 \leq i \leq \theta \}, \\
K_iK_i^{-1}&=L_iL_i^{-1}=1,  & \Et_k\Ft_\ell&=\Ft_\ell\Et_k,
\\ K_i\Et_k&=\chi(\alpha_i,\beta_k)\, \Et_k K_i,  & L_i \Et_k&=\chi(\beta_k,-\alpha_i)\,\Et_k L_i,
\\ K_i\Ft_k&=\chi(-\alpha_i,\beta_k)\,\Ft_kK_i, & L_i\Ft_k&=\chi(\beta_k,\alpha_i)\,\Ft_kL_i,
\\ \Et_k \Et_\ell&=\chi(\beta_k,\beta_\ell)\, \Et_\ell \Et_k, \quad k<\ell, & \Et_k^{N_k}&=0,\quad \beta_k\notin\cO(V),
\\ \Ft_k \Ft_\ell&=\chi(\beta_\ell,\beta_k)\, \Ft_\ell \Ft_k, \quad k<\ell, & \Ft_k^{N_k}&=0,\quad \beta_k\notin\cO(V).
\end{align*}
\end{corollary}
\pf The proof is analogous to the previous case if we check that $E_{\beta_k} F_{\beta_\ell}= F_{\beta_\ell}E_{\beta_k}$ in $\gr U(V)$; this relation follows by \eqref{eq:Fbeta Egamma}.
\epf

\medskip

\begin{theorem}\label{thm:noetheriano}
The algebras $U^{+}(V)$, $U^{\ge0}(V)$, $U(V)$ are Noetherian.
\end{theorem}
\pf
It suffices to prove that the graded algebras are Noetherian. As
\begin{itemize}
\item $\gr U^{+}(V)$ is the quotient of a quantum affine space on $M$ generators by powers of some generators,
\item $\gr U^{\ge0}(V)$ is the localization (we add the inverses of the $\theta$ generators $K_i$) of the quotient of a quantum affine space on $M+\theta$ generators by powers of some generators, and
\item $\gr U(V)$ is the localization of the quotient of a quantum affine space on $2M+2\theta$ generators by powers of some generators,
\end{itemize}
the three algebras $\gr U^{+}(V)$, $\gr U^{\ge0}(V)$, $\gr U(V)$ are Noetherian.
\epf

Now we compute the Gelfand-Kirillov dimension of these algebras. We refer to \cite{KL} for the definition and properties.

\begin{theorem}\label{thm:GKdim finita}
The Gelfand-Kirillov dimension of $U^{+}(V)$, $U^{\ge0}(V)$, $U(V)$ are, respectively, $|\cO(V)|$, $|\cO(V)|+\theta$, $2|\cO(V)|+2\theta$.
\end{theorem}
\pf
By \cite[Proposition 6.6]{KL}, $\GKdim U^{+}(V)=\GKdim \gr U^{+}(V)$.
The subalgebra $S^+(V)$ of $\gr U^{+}(V)$ generated by $\Et_k$, with $\beta_k\in\cO(V)$, is a quantum affine space in
$|\cO(V)|$ generators and $\gr U^{+}(V)$ is a free $S^+(V)$-module of rank $\prod\limits_{k:\beta_k\notin \cO(V)} N_k$, so $ \GKdim \gr U^{+}(V)=\GKdim S^+(V)= |\cO(V)|$.
For $U^{\ge 0}(V)$, let $S^{\ge0}(V)$ be the subalgebra of $\gr U^{\ge0}(V)$  generated by $K_i$, $i\in\I$, $\Et_k$ if $\beta_k\in\cO(V)$; it is the localization of a quantum affine space with $|\cO(V)|+\theta$ generators, and $\gr U^{\ge 0}(V)$ is a free $S^{\ge 0}(V)$-module of rank $\prod\limits_{k:\beta_k\notin \cO(V)} N_k$, so $ \GKdim U^{\ge 0}(V)=\GKdim \gr U^{\ge 0}(V)=\GKdim S^{\ge 0}(V)= |\cO(V)|+\theta$.
The proof for $U(V)$ follows analogously.
\epf

\section{Power root vectors on distinguished pre-Nichols algebras}\label{section:subalgebra Z}

Now we study the subalgebra $Z(V)$ of $U(V)$ generated by $E_\alpha^{N_\alpha}$, $F_\alpha^{N_\alpha}$, $K_\alpha^{N_\alpha}$, $L_\alpha^{N_\alpha}$, $\alpha\in\cO(V)$. First we describe the product of these elements. Then we give a general formula for the composition of the coproduct with $T_i^{V}$ on the whole $U^+_{-i}(V)$, and finally apply this formula to show that $Z(V)$ is a Hopf subalgebra. In what follows
\begin{itemize}
  \item $Z^+(V)$ and $Z^-(V)$ are the subalgebras generated by $E_\alpha^{N_\alpha}$, respectively $F_\alpha^{N_\alpha}$, $\alpha\in\cO(V)$.
  \item $\Gamma(V)$ is the subgroup of $\Z^{2\theta}$ generated by $K_\alpha^{N_\alpha}$, $L_\alpha^{N_\alpha}$, $\alpha\in\cO(V)$.
  \item $S(V)$ is a set of representatives of $\Z^{2\theta}/\Gamma(V)$.
\end{itemize}
The multiplication gives an isomorphism of $\zt$-graded vector spaces $Z^+(V)\ot \ku \Gamma(V)\ot Z^-(V)\simeq Z(V)$.

\subsection{Algebra structure of $Z(V)$}

Power root vectors are central elements of quantized enveloping algebras $U_q(\g)$, $q$ a root of unity \cite[Chapter 19]{DP-quantum}. In the general case they $q$-commute with the elements of $U(V)$.

\begin{proposition}\label{prop:power roots q-commute}
Let $\beta\in\cO(V)$.

\noindent\vi For each $X\in U^+(V)$ homogeneous of degree $\gamma$,
\begin{align*}
E_\beta^{N_\beta} X&= \chi(N_\beta\beta,\gamma)\,  XE_\beta^{N_\beta}, & F_\beta^{N_\beta} X&= XF_\beta^{N_\beta}.
\end{align*}

\noindent\vii For each $Y\in U^-(V)$ homogeneous of degree $\gamma$,
\begin{align*}
E_\beta^{N_\beta} Y&= YE_\beta^{N_\beta}, & F_\beta^{N_\beta} Y&= \chi(\gamma,N_\beta\beta)\, YF_\beta^{N_\beta}.
\end{align*}
\end{proposition}
\pf
If $\beta=\alpha_i$, then $E_i^{N_i} F_j=F_jE_i^{N_i}$ for $j\neq i$, and $E_i^{N_i} F_i=F_iE_i^{N_i}$ by \eqref{eq:simple power root contra Fp}, since $(N_i)_{q_{ii}}=0$. Therefore $E_i^{N_i} Y= YE_i^{N_i}$ for all $Y\in U^-(V)$.

For each $j\neq i$, $0=(\ad_c E_i)^{N_i} E_j= E_i^{N_i} E_j-q_{ij}^{N_i} E_j E_i^{N_i}$ by \eqref{eq:expresion adjunta}, since $N_i\geq 1-a_{ij}$.
Then $E_i^{N_i} X= \chi(N_i\alpha_i,\gamma) XE_i^{N_i}$ for all $X\in U^+(V)_\gamma$, $\gamma\in\N_0$.

If $\beta=s_{i_1}\cdots s_{i_{k-1}}(\alpha_{i_k})$, then $E_\beta=T_{i_1}\cdots T_{i_{k-1}}(\Eb_{i_k})$. Let $v=s_{i_{k-1}}\dots s_{i_1}$. If $j$ is such that $v(\alpha_j)\in\N_0^\theta$, then $T_{i_{k-1}}^-\dots T_{i_1}^-(E_j)\in U^+(V')_{v(\alpha_j)}$ by \cite[Proposition 6.15]{H-isom}, so we can use the previous case:

\begin{align*}
E_\beta^{N_\beta} E_j &= T_{i_1}\cdots T_{i_{k-1}}\left( v^* \chi(N_{i_k}\alpha_{i_k},v(\alpha_j)) \Eb_{i_k}^{N_{i_k}} T_{i_{k-1}}^-\dots T_{i_1}^-(E_j) \right) \\
&= \chi(N_\beta\beta,\alpha_j)  E_j E_\beta^{N_\beta}.
\end{align*}
As also $T_{i_{k-1}}^-\dots T_{i_1}^-(F_j)\in U^-(V')_{-v(\alpha_j)}$, we have that
\begin{align*}
E_\beta^{N_\beta} F_j &= T_{i_1}\cdots T_{i_{k-1}}\left(\Eb_{i_k}^{N_{i_k}} T_{i_{k-1}}^-\dots T_{i_1}^-(F_j) \right)     \\
&= T_{i_1}\cdots T_{i_{k-1}}\left(T_{i_{k-1}}^-\dots T_{i_1}^-(F_j) \Eb_{i_k}^{N_{i_k}}\right) = F_j E_\beta^{N_\beta}.
\end{align*}
If $j$ satisfies $v(\alpha_j)\in-\N_0^\theta$, then $s_{i_{t-1}}\cdots s_{i_1}(\alpha_j)\in\N_0^\theta$, $s_{i_t}\cdots s_{i_1}(\alpha_j)\in-\N_0^\theta$ for $t\geq 0$ minimal. Therefore $s_{i_{t-1}}\cdots s_{i_1}(\alpha_j)=\alpha_{i_t}$, so $T_{i_{t-1}}^-\dots T_{i_1}^-(E_j)=c\,E_{i_t}$ for some $c\in\ku^\times$, and
$$ T_{i_{k-1}}^-\dots T_{i_1}^-(E_j)= T_{i_{k-1}}^-\dots T_{i_{t+1}}^-(c\, K_t^{-1}F_t) \in T_{i_{k-1}}^-\dots T_{i_{t+1}}^-(K_t^{-1})U^-(V'), $$
because $s_{i_{k-1}}\dots s_{i_{t+1}}(-\alpha_t)=w(\alpha_j)\in-\N_0^\theta$, so
\begin{align*}
E_\beta^{N_\beta} & E_j = T_{i_1}\cdots T_{i_{k-1}}\left(\Eb_{i_k}^{N_{i_k}} T_{i_{k-1}}^-\dots T_{i_1}^-(E_j) \right) \\
&= c\, T_{i_1}\cdots T_{i_{k-1}}\left(\Eb_{i_k}^{N_{i_k}} T_{i_{k-1}}^-\dots T_{i_t}^-(cE_t) \right) \\
&= c\, T_{i_1}\cdots T_{i_{k-1}}\left(\Eb_{i_k}^{N_{i_k}} T_{i_{k-1}}^-\dots T_{i_{t+1}}^-(K_t^{-1})T_{i_{k-1}}^-\dots T_{i_{t+1}}^-(F_t) \right) \\
&= c\, T_{i_1}\cdots T_{i_{k-1}}\big( w^* \chi(N_{i_k}\alpha_{i_k},w(\alpha_j)) T_{i_{k-1}}^-\dots T_{i_{t+1}}^-(K_t^{-1}) \\
& \qquad \Eb_{i_k}^{N_{i_k}} T_{i_{k-1}}^-\dots T_{i_{t+1}}^-(F_t) \big) \\
&= c\, T_{i_1}\cdots T_{i_{k-1}}\big( w^* \chi(N_{i_k}\alpha_{i_k},w(\alpha_j)) T_{i_{k-1}}^-\dots T_{i_{t+1}}^-(K_t^{-1}) \\
& \qquad T_{i_{k-1}}^-\dots T_{i_{t+1}}^-(F_t) \Eb_{i_k}^{N_{i_k}} \big) = \chi(N_\beta\beta,\alpha_j)  E_j E_\beta^{N_\beta}.
\end{align*}
Analogously, $T_{i_{t-1}}^-\dots T_{i_1}^-(E_j)=c'\,E_{i_t}$ for some $c'\in\ku^\times$, so
$$ T_{i_{k-1}}^-\dots T_{i_1}^-(F_j)= T_{i_{k-1}}^-\dots T_{i_{t+1}}^-(c'\, E_t L_t^{-1}) \in U^+(V')T_{i_{k-1}}^-\dots T_{i_{t+1}}^-(L_t^{-1}), $$
because $s_{i_{k-1}}\dots s_{i_{t+1}}(-\alpha_t)=w(\alpha_j)\in-\N_0^\theta$, and

\begin{align*}
E_\beta^{N_\beta} & E_j = T_{i_1}\cdots T_{i_{k-1}}\left(\Eb_{i_k}^{N_{i_k}} T_{i_{k-1}}^-\dots T_{i_1}^-(E_j) \right) \\
& = T_{i_1}\cdots T_{i_{k-1}}\left(\Eb_{i_k}^{N_{i_k}} T_{i_{k-1}}^-\dots T_{i_t}^-(cE_t) \right) \\
&= c'\, T_{i_1}\cdots T_{i_{k-1}}\left(\Eb_{i_k}^{N_{i_k}} T_{i_{k-1}}^-\dots T_{i_{t+1}}^-(E_t)T_{i_{k-1}}^-\dots T_{i_{t+1}}^-(L_t^{-1}) \right) \\
&= c'\, T_{i_1}\cdots T_{i_{k-1}}\Big( w^* \chi(N_{i_k}\alpha_{i_k},w(\alpha_j))  T_{i_{k-1}}^-\dots T_{i_{t+1}}^-(E_t) \\
&\qquad \Eb_{i_k}^{N_{i_k}} T_{i_{k-1}}^-\dots T_{i_{t+1}}^-(L_t^{-1}) \Big) \\
&= c'\, T_{i_1}\cdots T_{i_{k-1}}\left( T_{i_{k-1}}^-\dots T_{i_{t+1}}^-(E_t) T_{i_{k-1}}^-\dots T_{i_{t+1}}^-(L_t^{-1}) \Eb_{i_k}^{N_{i_k}} \right) \\
&= \chi(N_\beta\beta,\alpha_j)  E_j E_\beta^{N_\beta},
\end{align*}
which completes the proof for $E_\beta^{N_\beta}$. The proof for $F_\beta^{N_\beta}$ is analogous.
\epf

\begin{corollary}\label{coro:derivadas PRV}
For all $\alpha\in\cO(V)$ and all $j\in\I$, $ \partial_j^K(E_\alpha^{N_\alpha})=\partial_j^L(E_\alpha^{N_\alpha})=0.$
\end{corollary}
\pf Let $i=i_1$.
The proof for $\alpha=\alpha_{i}$ is direct. If $\alpha\neq\alpha_i$, then set $\beta=s_i(\alpha)\in\cO(\rho_i(V))$, so  $E_\alpha=T_i(E_\beta)$, $N_\alpha=N_\beta$. By Proposition \ref{prop:power roots q-commute},
\begin{align*}
\partial_i^K(E_\alpha^{N_\alpha})K_i- &L_i^{-1}\partial_i^L(E_\alpha^{N_\alpha})= E_\alpha^{N_\alpha} F_i-F_iE_\alpha^{N_\alpha} \\
&= T_i\left( E_\beta^{N_\beta} E_iL_i^{-1} - E_iL_i^{-1}E_\beta^{N_\beta} \right) \\
&= T_i\left( s_i^*\chi( N_\beta\beta,\alpha_i) E_iE_\beta^{N_\beta} L_i^{-1} - E_iL_i^{-1}E_\beta^{N_\beta} \right)=0.
\end{align*}
Then $\partial_i^K(E_\alpha^{N_\alpha})=\partial_i^L(E_\alpha^{N_\alpha})=0$. For $j\neq i$,
\begin{align*}
\partial_j^K(E_\alpha^{N_\alpha})K_j- &L_j^{-1}\partial_j^L(E_\alpha^{N_\alpha})= E_\alpha^{N_\alpha} F_j-F_jE_\alpha^{N_\alpha} \\
&= T_{i}\left( E_\beta^{N_\beta} F_{i}^{-a_{ij}^{V}} F_j - F_{i}^{-a_{ij}^{V}} F_jE_\beta^{N_\beta} \right)=0,
\end{align*}
so $\partial_j^K(E_\alpha^{N_\alpha})=\partial_j^L(E_\alpha^{N_\alpha})=0$.
\epf

We deduce the freeness of $U(V)$ as $Z(V)$-module.

\begin{theorem}\label{thm:root vectors almost central, freeness}
The set
\begin{align*}
& \left\{ \left(\prod_{\beta\in\cO(V)}E_{\beta}^{k_\beta N_\beta}\right) \gamma  \left(\prod_{\beta\in\cO(V)}F_{\beta}^{l_\beta N_\beta}\right): k_\beta,l_\beta\in\N_0,\gamma\in\Gamma(V)  \right\},
\end{align*}
is a basis of $Z(V)$, where the order on the set $\cO(V)$ is arbitrary.

Moreover $U(V)$ is a free $Z(V)$-module with basis:
\begin{align*}
& \left\{ E_{\beta_M}^{a_M}E_{\beta_{M-1}}^{a_{M-1}} \cdots E_{\beta_1}^{a_1} s F_{\beta_M}^{b_M}F_{\beta_{M-1}}^{b_{M-1}} \cdots
F_{\beta_1}^{b_1}\, |  0\leq a_k,b_k < N_k , s\in S(V) \right\}.
\end{align*}
\end{theorem}
\pf
By Proposition \ref{prop:power roots q-commute}, the generators of the subalgebra $Z(V)$ $q$-commute with all the elements of $U(V)$.
Then apply Theorem \ref{thm: PBW bases preNichols}.
\epf

\begin{remark}
The algebra $Z(V)$ is not necessarily central since we can have $\chi(N_\beta\beta,\alpha_j)\neq1$ for some $\beta\in\cO(V)$, $1\le j\le \theta$.
\end{remark}

\begin{lemma}\label{lema:caracter sim trivial en Nbeta beta}
Let $\beta\in\cO(V)$. For every $\alpha\in\Z^\theta$, $\chi(N_\beta\beta,\alpha)\chi(\alpha,N_\beta\beta)=1$.
\end{lemma}
\pf
If $\beta=\alpha_i$, then $\chi(N_i\alpha_i,\alpha_i)\chi(\alpha_i,N_i\alpha_i)=q_{ii}^{2 N_i}=1$ and for $j\neq i$,
$$\chi(N_i\alpha_i,\alpha_j)\chi(\alpha_j,N_i\alpha_i)=(q_{ij}q_{ji})^{N_i}=q_{ii}^{c_{ij}^{V} N_i}=1,$$
so $\chi(N_i\alpha_i,\alpha)\chi(\alpha,N_i\alpha_i)=1$ for every $\alpha\in\Z^\theta$.
If $\beta\neq\alpha_i$, then $\beta=w(\alpha_i)$ for some $w\in\Hom(\cW,V)$ and $\alpha_i\in\cO(w^{-1}(V))$, so $N_\beta=N_i$. Therefore
\begin{align*}
\chi(N_\beta\beta,\alpha)\chi & (\alpha,N_\beta\beta)=\chi(N_i w(\alpha_i),\alpha)\chi(\alpha,N_i w(\alpha_i))\\
&=(w^{-1})^*\chi(N_i \alpha_i,w^{-1}(\alpha)) (w^{-1})^*\chi(w^{-1}(\alpha),N_i \alpha_i)=1,
\end{align*}
for all $\alpha\in\Z^\theta$, by applying the previous step to $\alpha_i\in\cO(w^{-1}(V))$.
\epf

\subsection{On the coproduct of $U(V)$}
We want to factorize the composition $\underline{\Delta}\circ T_i: U_{-i}^+(V)\to U^+(V)\ot U_{+i}^+(V)$ as in \cite{HS-coideal subalg}. Similar formula was introduced in \cite[Proposition 5.3.4]{L-libro} for quantized enveloping algebras. Such factorization is interpreted as equivalences between the corresponding categories of Yetter-Drinfeld modules \cite{BLS,HS-general}. One of the factors is $T_i\ot T_i$ but viewed as an algebra map between \emph{braided} structures.
\begin{itemize}
  \item $U_{-i}^+(V)\underline{\ot} U_{-i}^+(V)$ denotes the tensor product $U_{-i}^+(V)\ot U_{-i}^+(V)$ with the algebra structure in $^{\cB_{i}^{\bop}\#\ku\Z^\theta}_{\cB_{i}^{\bop}\#\ku\Z^\theta}\mathcal{YD}$, see Remark \ref{rem:algebra de Hopf YD para U-i}, and
  \item $U_{+i}^+(V)\underline{\ot} U_{+i}^+(V)$ denotes the tensor product $U_{+i}^+(V)\ot U_{+i}^+(V)$ with the algebra structure in $^{\cB_{i}\#\ku\Z^\theta}_{\cB_{i}\#\ku\Z^\theta}\mathcal{YD}$, see Remark \ref{rem:algebra de Hopf YD para U+i}.
\end{itemize}

\begin{lemma}\label{lemma: Ti ot Ti es de algebras}
$T_i\ot T_i: U_{-i}^+(V)\underline{\ot} U_{-i}^+(V)\to U_{+i}^+(V)\underline{\ot} U_{+i}^+(V)$ is an algebra map.
\end{lemma}
\pf
Let $x,y,w,z\in U_{-i}^+(V)$, where $y$, $w$ are homogeneous of degrees $\beta,\gamma\in\N_0^\theta$, respectively. As $(\ad_c E_i)^{N_i}\equiv 0$,

\begin{align*}
T_i & \ot T_i  (x\ot y) \cdot \, T_i\ot T_i(w\ot z)= T_i(x) \left( T_i(y)_{(-1)}\rightharpoonup T_i(w) \right) \underline{\ot} T_i(y)_{(0)}  T_i(z) \\
& \overset{\eqref{eq:accion-coaccion a izquierda}}= \sum_{n=0}^{N_i-1} s_i^*\chi(s_i(\beta)-n\alpha_i,s_i(\gamma)) \, T_i(x) (\ad_{\cb} \Eb_i)^n T_i(w)\underline{\ot} \, (T_i(y) \triangleleft \Ef_i^{(n)}) T_i(z) \\
& \overset{\eqref{eq:producto coproducto dual},\eqref{eq:Ti con delta K}}= \sum_{n=0}^{N_i-1} \frac{(-1)^nq_{ii}^{\frac{n(n-1)}{2}}\chi(\beta+n\alpha_i,\gamma)}{(n)_{q_{ii}}! \chi(\alpha_i,\gamma)^{n}}
\, T_i(x) T_i((\partial_i^K)^n w) \\
& \qquad \qquad \qquad \underline{\ot} \left((\partial_i^L)^n T_i(y) \right) T_i(z) \\
& \overset{\eqref{eq:producto coproducto dual}}= \sum_{n=0}^{N_i-1} (-1)^nq_{ii}^{\frac{n(n-1)}{2}}\chi(\beta,\gamma)
\, T_i(x) T_i(\Ef_i^{(n)} \triangleright w)\underline{\ot} \left((\partial_i^L)^n T_i(y) \right) T_i(z) \\
& \overset{\eqref{eq:Ti con delta L}}= \sum_{n=0}^{N_i-1} \frac{\chi(\beta,\gamma)}{\chi(\beta,\alpha_i)^{n}}
\, T_i(x) T_i(\Ef_i^{(n)} \triangleright w)\underline{\ot} T_i \left(E_i^n \rightharpoonup ' y  \right) T_i(z) \\
& = \sum_{n=0}^{N_i-1} \chi(\beta,\gamma-n\alpha_i)
\, T_i \left( x \left(\Ef_i^{(n)} \triangleright w \right) \right)\underline{\ot} T_i \left( \left(E_i^n \rightharpoonup ' y \right) z \right) \\
& = T_i\ot T_i \left( (x\ot y) \cdot (w\ot z) \right),
\end{align*}
so the proof is complete.
\epf

Let $Y\in U^+_{+i}(V)_\beta$, $\beta\in\N_0^\theta$. As $Y \triangleleft\Ef_i^{(n)}\in U(V)_{\beta-n\alpha_i}$, there exists $k\in\N$ such that $Y \triangleleft\Ef_i^{(n)} =0$ for all $n\geq k$. Then there exists a well-defined map $\Rf_i:U_{+i}^+(V)\ot U_{+i}^+(V)\to U^+(V)\ot U_{+i}^+(V)$ such that
\begin{align}\label{eq:def Ri}
\Rf_i^{V}(x\otimes y)&= \sum_{k\geq0}  x E_i^k\ot y \triangleleft\Ef_i^{(k)}, & x,y &\in U^+_{+i}(V).
\end{align}
Compare with \cite[(3.4)]{HS-coideal subalg}.

\begin{lemma}\label{lemma:Rfi es de algebras}
$\Rf_i^{V}:U_{+i}^+(V)\underline{\ot} U_{+i}^+(V) \to U^+(V)\ot U_{+i}^+(V)$ is an algebra map.
\end{lemma}

Here $U^+(V)\ot U_{+i}^+(V)$ is a subalgebra of $U^+(V)\ot U^+(V)$, where $U^+(V)\ot U^+(V)$ has the algebra structure such that
$ (x\ot y)(x'\ot y)= \chi(\beta,\gamma) xx'\ot yy$
for $x,x',y,y'\in U^+(V)$, $y$, $x'$ homogeneous of degree $\beta$, $\gamma$, respectively.

\pf
Set $x,y,w,z\in U_{+i}^+(V)$, $y$, $w$ homogeneous of degrees $\beta,\gamma\in\N_0^\theta$. Then
\begin{align}\label{eq:formula adjunta}
E_i^n \, w
&= \sum_{r=0}^n \binom{n}{r}_{q_{ii}}\chi(\alpha_i,\gamma)^{n-r} (\ad_c E_i)^r w \, E_i^{n-r}.
\end{align}
for all $n\in\N$, so

\begin{align*}
\Rf_i^{V} &(x\ot y) \cdot \Rf_i^{V}(w\ot z)= \\
&= \sum_{j,k\geq 0} \chi(\beta-j\alpha_i,\gamma+k\alpha_i) \, xE_i^j w E_i^k \ot (y\triangleleft\Ef_i^{(j)}) (z\triangleleft\Ef_i^{(k)}) \\
& \overset{\eqref{eq:formula adjunta}} = \sum_{j,k\geq 0}
\, x  \left( \sum_{r=0}^{\min \{N_i-1,j\} } \binom{j}{r}_{q_{ii}}
\frac{\chi(\beta,\gamma)\chi(\beta,\alpha_i)^k}{\chi(\alpha_i,\gamma)^rq_{ii}^{jk} }
 (\ad_c E_i)^r w E_i^{k+j-r} \right)
\\ &\qquad \qquad\qquad \ot (y\triangleleft\Ef_i^{(j)}) (z\triangleleft\Ef_i^{(k)})
\\
& \overset{\eqref{eq:producto coproducto dual}} = \sum_{j,k\geq 0} \sum_{r=0}^{\min \{N_i-1,j\} }
\frac{\chi(\beta,\gamma+k\alpha_i)}{\chi(\alpha_i,\gamma)^r q_{ii}^{jk} }
x(\ad_c E_i)^r w E_i^{k+j-r} \\
&\qquad \qquad\qquad \ot ((y\triangleleft\Ef_i^{(r)}) \triangleleft\Ef_i^{(j-r)}) (z\triangleleft\Ef_i^{(k)}).
\end{align*}
On the other hand,
\begin{align*}
\Rf_i^{V} &\left( (x\ot y) \cdot (w\ot z)\right) =\\
&\overset{\eqref{eq:accion-coaccion a izquierda}}=  \Rf_i^{V} \left( \sum_{r=0}^{N_i} \chi(\beta-r\alpha_i,\gamma) \,  x(\ad_c E_i)^r w \ot (y \triangleleft \Ef_i^{(r)})z \right) \\
& \overset{\eqref{eq:def Ri}} = \sum_{r=0}^{N_i}\sum_{t\geq 0} \chi(\beta-r\alpha_i,\gamma) \,  x(\ad_c E_i)^r w E_i^t \ot \left((y \triangleleft \Ef_i^{(r)})z \right)\triangleleft \Ef_i^{(t)}
\\ & \overset{\eqref{eq:Efi como skew derivation}} = \sum_{r=0}^{N_i}\sum_{k,l\geq 0} \chi(\beta-r\alpha_i,\gamma)\chi(\beta-(r+l)\alpha_i,k\alpha_i) \,  x(\ad_c E_i)^r w E_i^{k+l}
\\ & \qquad \qquad \qquad \ot ((y \triangleleft \Ef_i^{(r)}) \triangleleft \Ef_i^{(l)} )z \triangleleft \Ef_i^{(k)}
\\ & = \sum_{r=0}^{N_i}\sum_{k,l\geq 0} \frac{\chi(\beta,\gamma+k\alpha_i)}{\chi(\alpha_i,\gamma)^r q_{ii}^{(r+l)k}} \,  x(\ad_c E_i)^r w E_i^{k+l} \ot ((y \triangleleft \Ef_i^{(r)}) \triangleleft \Ef_i^{(l)} )z \triangleleft \Ef_i^{(k)},
\end{align*}
so $\Rf_i^{V} (x\ot y) \cdot \Rf_i^{V}(w\ot z)=\Rf_i^{V} \left( (x\ot y) \cdot (w\ot z)\right)$.
\epf

Recall the map $\Delta_i:U_{+i}^+(V) \to U_{+i}^+(V)\underline{\ot} U_{+i}^+(V)$ introduced in Remark \ref{rem:algebra de Hopf YD para U-i}. Now we prove the factorization of $\underline{\Delta}\circ T_i$.

\begin{theorem}\label{thm:factorization delta Ti}
For every $E\in U_{-i}^+(V)$,
\begin{equation}\label{eq:factorization delta Ti}
\underline{\Delta}\circ T_i (E)= \Rf_i^{\rho_i(V)} \circ (T_i\ot T_i)\circ \Delta_i(E).
\end{equation}
\end{theorem}
\pf
It is enough to prove the formula for $E=E_{j,n}^-$, $j\neq i$, $0\leq n\leq -c_{ij}^{V}$, since they generate $U_{-i}^+(V)$ as an algebra and both $\underline{\Delta}\circ T_i$, $\Rf_i^{\rho_i(V)} \circ (T_i\ot T_i)\circ \Delta_i$ are algebra maps by Lemmas \ref{lemma: Ti ot Ti es de algebras} and \ref{lemma:Rfi es de algebras}. Let
$$ \kappa_{i,j;n}:= \underline{q}_{ii}^{\,\, n} \left(\prod_{t=-c_{ij}^{V}-n}^{-c_{ij}^{V}-1} (t+1)_{\underline{q}_{ii}}(1-\underline{q}_{ii}^{\,\, t}\underline{q}_{ij}\underline{q}_{ji}) \right).$$
By \eqref{eq:Ti(E-)=E+} and \eqref{eq:coproducto Ei+},
\begin{align*}
\kappa_{i,j;n}^{-1} &  \underline{\Delta}\circ T_i(E_{j,n}^-) =  \underline{\Delta}(\Eb_{j,-c_{ij}^{V}-n}^+) =  \Eb_{j,-c_{ij}^{V}-n}^+ \ot 1 \\
+& \sum_{s=0}^{-c_{ij}^{V}-n-1} \binom{-c_{ij}^{V}-n}{s}_{\underline{q}_{ii}}\left( \prod_{r=1}^s (1-\underline{q}_{ii}^{-c_{ij}^{V}-n-r}\underline{q}_{ij}\underline{q}_{ji})  \right) \Eb_i^{s} \otimes \Eb_{j,-c_{ij}^{V}-n-s}^+
\end{align*}
On the other hand, by \eqref{eq:formula delta sub i}, \eqref{eq:Ti(E-)=E+}, \eqref{eq:def Ri} and \eqref{eq:producto coproducto dual},

\begin{align*}
\kappa_{i,j;n}^{-1} \, \Rf_i^{\rho_i(V)} & \circ (T_i\ot T_i)\circ \Delta_i(E_{j,n}^-)=  \Rf_i^{\rho_i(V)} \left( \Eb_{j,-c_{ij}^{V}-n}^+\ot 1 + 1 \ot \Eb_{j,-c_{ij}^{V}-n}^+ \right) \\
&= \Eb_{j,-c_{ij}^{V}-n}^+\ot 1 + \sum_{s=0}^{-c_{ij}^{V}-n-1} \frac{1}{s_{\underline{q}_{ii}}!} \Eb_i^{s} \otimes (\partial_i^L)^s (\Eb_{j,-c_{ij}^{V}-n}^+),
\end{align*}
so $\underline{\Delta} \circ T_i(E_{j,n}^-) = \Rf_i^{\rho_i(V)} \circ (T_i\ot T_i)\circ \Delta_i(E_{j,n}^-)$.
\epf

\subsection{Coproduct structure of $Z(V)$}

Next we prove that $Z(V)$ is a Hopf subalgebra of $U(V)$. We start by describing the left hand side of the coproduct $\underline{\Delta}$ of $Z^+(V)$. Let $v=s_{i_1}\cdots s_{i_k}\in\Hom(\cW,V)$ be an element of length $k$, $\beta_j=s_{i_1}\cdots s_{i_{j-1}}(\alpha_{i_j})\in\Delta_+^{V}$, and as above $E_{\beta_j}=T_{i_1}\cdots T_{i_{j-1}}(E_{i_j})$.

\begin{proposition}\label{prop:lado izquierdo Delta Z(chi)}
If $\beta_k\in\cO(V)$, then there exist $X(n_1,\dots,n_{k-1})\in U^+(V)$ such that
\begin{align}\label{eq:formula lado izquierdo Delta Z(chi)}
\underline{\Delta}(E_{\beta_k}^{N_{\beta_k}}) &= E_{\beta_k}^{N_{\beta_k}}\ot 1+1\ot E_{\beta_k}^{N_{\beta_k}} \\
&\quad + \sum_{n_i\in \N_0} E_{\beta_{k-1}}^{n_{k-1}N_{\beta_{k-1}}} \dots E_{\beta_1}^{n_1N_{\beta_1}} \otimes X(n_1,\dots,n_{k-1}). \notag
\end{align}
\end{proposition}
\pf
The proof is by induction on $k$. If $k=1$, then $E_{\beta_1}=E_{i_1}$, and $E_{i_1}^{N_{i_1}}$ is primitive.
Assume that \eqref{eq:formula lado izquierdo Delta Z(chi)} holds for every $v'$ of length less than $k$. Let $v=s_{i_1}^{V} v'$, $i=i_1$, $\beta=\beta_k$, $\gamma_j:=s_{i_2}\cdots s_{i_{j-1}}(\alpha_{i_j})\in\Delta_+^{\rho_i(V)}$, $j=2,\cdots,k$, $\gamma=\gamma_k$, so $\beta_j=s_i(\gamma_j)$, $E_{\beta_j}=T_i(E_{\gamma_j})$, $N_\beta=N_\gamma$. By inductive hypothesis,
\begin{align*}
\underline{\Delta}(E_{\gamma}^{N_{\gamma}}) &= E_{\gamma}^{N_{\gamma}}\ot 1+1\ot E_{\gamma}^{N_{\gamma}} \\
&\qquad + \sum_{n_i\in \N} E_{\gamma_{k-1}}^{n_{k-1}N_{\gamma_{k-1}}} \dots  E_{\gamma_2}^{n_2N_{\gamma_2}}\otimes Y(n_1,\dots,n_{k-1}),
\end{align*}
for some $Y(n_1,\dots,n_{k-1})\in U^+(\rho_i(V))$. As $\partial_i^L(E_{\beta_j}^{N_{\beta_j}})=0$ by Corollary \ref{coro:derivadas PRV}, we have that $E_i\rightharpoonup' E_{\gamma_j}^{N_{\gamma_j}}=0$ by \eqref{eq:Ti con delta L}. By Lemma \ref{lema:caracter sim trivial en Nbeta beta},

\begin{align*}
0 &= \underline{\Delta}(E_i\rightharpoonup' E_{\gamma}^{N_{\gamma}}) = \underline{\Delta}(E_{\gamma}^{N_{\gamma}})(E_i\ot 1+1\ot E_i) \\
& \qquad -\chi(N_\gamma \gamma,\alpha_i) (E_i\ot 1+1\ot E_i) \underline{\Delta}(E_{\gamma}^{N_{\gamma}})  \\
&= \sum_{n_i\in \N} E_i\rightharpoonup' \left( E_{\gamma_{k-1}}^{n_{k-1}N_{\gamma_{k-1}}} \dots  E_{\gamma_2}^{n_2N_{\gamma_2}}\right) \otimes Y(n_2,\dots,n_{k-1}) \\
& \quad + E_{\gamma_{k-1}}^{n_{k-1}N_{\gamma_{k-1}}} \dots E_{\gamma_2}^{n_2N_{\gamma_2}} \otimes E_i\rightharpoonup' Y(n_2,\dots,n_{k-1}) \\
& \quad + \left(1- \prod_{j=2}^{k-1} \chi(N_{\gamma_j}\gamma_j,\alpha_i)^{n_j}\chi(\alpha_i,N_{\gamma_j}\gamma_j)^{n_j} \right)
\chi \left(N_{\gamma}\gamma-\sum_{j=2}^{k-1} N_{\gamma_j}\gamma_j,\alpha_i \right)  \\
& \qquad E_{\gamma_{k-1}}^{n_{k-1}N_{\gamma_{k-1}}} \dots  E_{\gamma_2}^{n_2N_{\gamma_2}} \otimes  Y(n_2,\dots,n_{k-1}) E_i \\
&= \sum_{n_i\in \N} E_{\gamma_{k-1}}^{n_{k-1}N_{\gamma_{k-1}}} \dots  E_{\gamma_2}^{n_2N_{\gamma_2}} \otimes E_i\rightharpoonup' Y(n_2,\dots,n_{k-1}),
\end{align*}
so $E_i\rightharpoonup' Y(n_2,\dots,n_{k-1})=0$ for all $n_i\in\N_0$. Moreover $(\underline{\Delta}\ot \id)\underline{\Delta}(E_{\gamma}^{N_{\gamma}})$ can be written as a sum of terms where the three elements of $U^+(V)^{\ot 3}$ are products of $E_{\gamma_j}^{n_jN_{\gamma_j}}$ and $Y(n_2,\dots,n_{j-1})$, $1\leq j\leq k$, so they are annihilated by $E_i\rightharpoonup'$. In particular the middle term have terms $E_i^n$ only for $n \in \N N_i$, so the remaining terms in $\iota_{i,V} \mathcal{S}\pi_{i,V}(x_{(2)})$ are of the form $E_i^{nN_i}$, $n\in\N_0$.
As $\Delta_i(E_{\gamma}^{N_{\gamma}})= (\id\ot m)(\id\ot\iota_{i,V} \mathcal{S}\pi_{i,V}\ot\id)(\underline{\Delta}\ot \id)\underline{\Delta}(E_{\gamma}^{N_{\gamma}})$, there exist $Z(n_2,\dots,n_{k-1})\in U^+(V)$ such that $E_i\rightharpoonup' Z(n_2,\dots,n_{k-1})=0$ and
\begin{align*}
\Delta_i(E_{\gamma}^{N_{\gamma}})& = E_{\gamma}^{N_{\gamma}}\ot 1+1\ot E_{\gamma}^{N_{\gamma}} \\
&\qquad + \sum_{n_i\in \N} E_{\gamma_{k-1}}^{n_{k-1}N_{\gamma_{k-1}}} \dots  E_{\gamma_2}^{n_2N_{\gamma_2}} \otimes  Z(n_2,\dots,n_{k-1}).
\end{align*}
By Theorem \ref{thm:factorization delta Ti},

\begin{align*}
\underline{\Delta}(& E_{\beta}^{N_{\beta}}) = \underline{\Delta}\circ T_i(E_{\gamma}^{N_{\gamma}})=\Rf_i^{V} \circ (T_i\ot T_i)\circ\Delta_i(E_{\gamma}^{N_{\gamma}}) \\
&= \Rf_i^{V} \left( E_{\beta}^{N_{\beta}}\ot 1+1\ot E_{\beta}^{N_{\beta}} \right) \\
&\qquad + \Rf_i^{V} \left( \sum_{n_i\in \N} E_{\beta_{k-1}}^{n_{k-1}N_{\beta_{k-1}}} \dots E_{\beta_2}^{n_2N_{\beta_2}} \otimes T_i(Z(n_2,\dots,n_{k-1})) \right).
\end{align*}
If $n$ is not a multiple of $N_i$, then
\begin{align*}
( T_i & (Z(n_2,\dots,n_{k-1})  ) )  \triangleleft\Ef_i^{(n)} \overset{\eqref{eq:producto coproducto dual}}= (n)_{q_{ii}}^{-1} \left( \partial_i^L \circ T_i(Z(n_2,\dots,n_{k-1})) \right)\triangleleft\Ef_i^{(n-1)} \\
& \overset{\eqref{eq:Ti con delta L}}= -\chi(\delta,\alpha_i)^{-1}q_{ii}^{-1}(n)_{q_{ii}}^{-1} \left(T_i(E_i\rightharpoonup'Z(n_2,\dots,n_{k-1})) \right)\triangleleft\Ef_i^{(n-1)} =0,
\end{align*}
where $\delta=N_\gamma \gamma- \sum_{j=2}^{k-1} n_jN_{\gamma_j}\gamma_j$ is the degree of $Z(n_2,\dots,n_{k-1})$. As also $\delta_i^L(E_{\beta}^{N_{\beta}})=0$, we have that $E_{\beta}^{N_{\beta}}\triangleleft\Ef_i^{(n)}=0$. Therefore

\begin{align*}
\underline{\Delta}(& E_{\beta}^{N_{\beta}}) = E_{\beta}^{N_{\beta}}  \ot 1 + \sum_{n\in\N_0} E_i^{n N_i}\ot E_{\beta}^{N_{\beta}} \triangleleft\Ef_i^{(nN_i)} \\
&  + \sum_{n_i\in \N}  \sum_{n\in\N_0} E_{\beta_{k-1}}^{n_{k-1}N_{\beta_{k-1}}} \dots E_{\beta_2}^{n_2N_{\beta_2}} E_i^{n N_i} \otimes T_i(Z(n_2,\dots,n_{k-1})) \triangleleft\Ef_i^{(nN_i)},
\end{align*}
which concludes the inductive step. \epf


\begin{theorem}\label{thm:Z+ es coinv del ker pi}
$Z^+(V)$ is a normal right coideal subalgebra of $U^+(V)$.

Moreover $Z^+(V)=^{co \pi_V|_{U^+(V)}}U^{+}(V).$
\end{theorem}
\pf
By Proposition \ref{prop:lado izquierdo Delta Z(chi)} $Z^+(V)$ is a right coideal subalgebra, and is normal by Theorem \ref{thm:root vectors almost central, freeness}. Then the right coideal $U^+(V) (Z^+(V))^+$ is a Hopf ideal, and $Z^+(V)=^{co \pi_V|_{U^+(V)}}U^{+}(V)$ by Proposition \ref{pro:masuoka} \vii, since $\u^{\ge 0}(V) = U^+(V)/ U^+(V) (Z^+(V))^+$.
\epf

We consider the following subalgebras of $U^+(V)$:
\begin{align*}
\cK^K(V) & =\bigcap_{i=1}^\theta \ker \partial_i^K, & \cK^L(V) & =\bigcap_{i=1}^\theta \ker \partial_i^L, & \cK(V)&= \cK^K(V) \cap\cK^L(V).
\end{align*}

\begin{remark}
$\cK^L(V)$ is a right coideal subalgebra since it is the intersection of right coideal subalgebras. Similarly $\cK^K(V)$ is a left coideal subalgebra.
\end{remark}

\begin{lemma}\label{lema:S entre ker derivaciones}
\vi For all $x\in U^+(V)$, $i\in\I$, $\partial_i^L(\mathcal{S}(x))= -\mathcal{S}(\partial_i^K(x))$.

\smallbreak
\noindent\vii The restriction of $\mathcal{S}$ gives bijective $\zt$-graded antialgebra maps
\begin{align*}
\ker \partial_i^K   & \overset{\sim}\longrightarrow \ker \partial_i^L, \quad i\in\I, &  \cK^K(V) & \longrightarrow \cK^L(V).
\end{align*}
\end{lemma}
\pf
\vi Let $x\in U^+(V)_n$, $n\in\N_0$. If $n=0$, then both sides are 0. Let $n\ge1$. As $\mathcal{S}$ is an anticoalgebra graded map and $\mathcal{S}(E_i)=-E_i$ for all $i\in\I$,
\begin{align*}
\sum_{i=1}^\theta \partial_i^L(\mathcal{S}(x)) \ot E_i &= \underline{\Delta}_{1,n-1}^{\cop}\circ\mathcal{S}(x)=(\mathcal{S}\ot\mathcal{S})\circ\underline{\Delta}_{n-1,1}(x) \\
&= -\sum_{i=1}^\theta   \mathcal{S}(\partial_i^K(x)) \ot E_i.
\end{align*}

\noindent\vii From the identity in \vi, $\mathcal{S}(\ker \partial_i^K)\subseteq \ker \partial_i^L$, $\mathcal{S}(\cK^K(V))\subseteq\cK^L(V)$. Also $\mathcal{S}$ is a $\zt$-graded, and bijective since $U^+(V)$ is connected.
\epf

We obtain a different characterization of
$Z^+(V)$.

\begin{theorem}\label{thm:Z+ es subalg hopf U+}
The subalgebra $Z^+(V)$ coincides with $\cK^L(V)$, $\cK^K(V)$. In particular $Z^+(V)$ is a braided Hopf subalgebra of $U^+(V)$.
\end{theorem}
\pf
By Corollary \ref{coro:derivadas PRV}, $Z^+(V) \subseteq \cK^L(V)$. If $x\in \cK^L(V)$ is an homogeneous element of degree $n\geq2$, then $\pi_V(x)$ is an homogeneous element of the same degree, annihilated by all the corresponding derivations $\partial_i^L$ of the Nichols algebra $\u^+(V)$. But the intersection of the kernels of $\partial_i^L$ in $\u^+(V)$ is $\ku 1$ by \cite[Proposition 2.4]{MS}, so $\pi_V(x)=0$. Then
$ U^+(V) N(\cK^L(V))^+ \subseteq \ker\pi_V$. As $\ker\pi_V=U^+(V) (Z^+(V))^+$ by Theorem \ref{thm:Z+ es coinv del ker pi}, both Hopf ideals coincide, so by Proposition \ref{pro:masuoka} $Z^+(V)=N(\cK^L(V))$. Therefore $Z^+(V)=\cK^L(V)$.

Also $Z^+(V) \subseteq \cK^K(V)$ by Corollary \ref{coro:derivadas PRV}, and by Lemma \ref{lema:S entre ker derivaciones}\vii the Hilbert series of $\cK^K(V)$ and $\cK^L(V)$ coincide. As $Z^+(V)=\cK^L(V)$, we have that $Z^+(V)=\cK^K(V)$. Therefore $Z^+(V)$ is a braided Hopf subalgebra of $U^+(V)$ since it is simultaneously a right and a left coideal subalgebra.
\epf

\begin{corollary}
The subalgebras $Z^\pm(V)$, $Z^0(V)$, $Z(V)$ do not depend on the choice of the element of maximal length.
\end{corollary}
\pf
By Theorem \ref{thm:Z+ es subalg hopf U+}, $Z^+(V)=\cK(V)$, and $\cK(V)$ does not depend on that element of maximal length; the proof for $Z^-(V)$ is analogous. As $Z^0(V)$ is the group algebra of the subgroup generated by $K_\beta^{N_\beta}$, $L_\beta^{N_\beta}$, $\beta\in\cO(V)$, the proof for $Z(V)$ follows since it is generated by $Z^\pm(V)$, $Z^0(V)$.
\epf

We now deduce the corresponding statement for $Z(V)$ and $U(V)$.

\begin{theorem}\label{thm:Z es subalg hopf U}
$Z(V)$ is a Hopf subalgebra of $U(V)$.
\end{theorem}
\pf
Let $\beta\in\cO(V)$. The $U^{+0}(V)$-coaction is $\lambda(E_\beta^{N_\gamma})=K_\gamma^{N_\gamma}\ot E_\gamma^{N_\gamma}$, $\gamma\in\cO(V)$, and $\underline{\Delta}(Z^+(V))\subset Z^+(V)\ot Z^+(V)$ by Theorem \ref{thm:Z+ es subalg hopf U+}, so $\Delta(Z^+(V))\subset Z^+(V)\Gamma(V)\ot Z^+(V)$. Analogously, $\Delta(Z^-(V))\subset Z^-(V)\ot Z^-(V)\Gamma(V)$.
\epf

\begin{remark}
Assume that $\chi$ is such that $\chi(N_\beta\beta,\alpha_j)$=1 for all $\beta\in\cO(V)$, $1\le j\le\theta$. Then $Z(V)$ is a central Hopf subalgebra and $\Gamma(V)$ acts trivially on $U^\pm(V)$, $\u^\pm(V)$ so the quotient
\begin{align*}
\widehat{\u}(V)=\u(V)/<K^{N_\beta\beta}-1,L^{N_\beta\beta}-1:\beta\in\cO(V)>
\end{align*}
is well defined. The triangular decomposition on $\u(V)$ induces a triangular decomposition on $\widehat{\u}(V)$, where $\widehat{\u}^\pm(V)$ is canonically identified with $\u^\pm(V)$, and the zero part is $\widehat{\u}^0(V)\simeq \ku(\Z^{2\theta}/\Gamma(V))$. The algebra $\widehat{\u}(V)$ can be seen as a quantum double of the bosonizations of $\u^\pm(V)$ with $\ku(\zt/<K^{N_\beta\beta}-1:\beta\in\cO(V)>)$. We have an extension of Hopf algebras \cite{AD}:
$$ \ku \rightarrow Z(V) \hookrightarrow U(V) \twoheadrightarrow \widehat{\u}(V) \rightarrow \ku, $$
which generalizes the case of quantum groups $U_q(\g)$ at roots of unity and the corresponding small quantum groups $\u_q(\g)$ obtained as a quotient.
\end{remark}

\begin{remark}
The graded dual $\mathcal{L}(V)=U^+(V)^*$ is a braided Hopf algebra containing a copy of $\u^+(V)^* \simeq \u^+(V)$; moreover is a finite module over this subalgebra. $\mathcal{L}(V)$ was called the \emph{Lusztig algebra} \cite{AAGTV} and resembles the $q$-divided power algebra \cite{L-libro}.
\end{remark}

\section{Some examples}\label{section:examples}

\subsection{Braidings of super type $A$}
We apply previous results to braidings of super type $A$, see \cite{AAY}. Fix $N\in\N$, $q\in\ku$ a root of unity of order $N\geq 3$. There exists a braiding $C(\theta,q;i_1,\ldots ,i_j)$ of super type $A_\theta$ for each ordered subset $1\leq i_1<\dots<i_k\leq\theta$. Its matrix $(q_{ij})_{1 \leq i,j \leq \theta}$ satisfies
\begin{itemize}
\smallbreak\item $q_{ij}q_{ji}=1$ if $1 \le i, j<\theta$, $|i-j|>1$,
\smallbreak\item if $i=i_\ell$ for some $1\leq\ell\le k$, then $q_{ii}=-1$, $q_{i-1,i}q_{i,i-1}q_{i+1,i}q_{i,i+1}=1$,
\smallbreak\item if $i\neq i_\ell$ for all $1\leq\ell\le k$, then $q_{ii}q_{i-1,i}q_{i,i-1}=q_{ii}q_{i+1,i}q_{i,i+1}=1$,
\end{itemize}
where $q=q_{11}^2 q_{12}q_{21}$. The notation is close to the simple chains in \cite{H-classif RS}.

The positive roots corresponding to this braiding are
\begin{align*}
\Delta_+^{V} &= \{ \alpha_{j,k}:1\le j\le k\le\theta\}, & \mbox{where }  \alpha_{j,k} & =\alpha_j+\alpha_{j+1} +\dots + \alpha_k.
\end{align*}
Set $t_k^{V}=s_1^{V} s_2\cdots s_k$,  $w=t_\theta^{V} t_{\theta-1}\dots t_1$, where we use for the $t$'s the same convention for concatenation as for the $s$'s. The expression of $w$ in $s_i$'s is reduced and this is the element of maximal length by Lemma \ref{Lemma:longitudHY} since
\begin{align}
t_\theta &t_{\theta-1}\dots t_{\theta-j+2} s_1\dots s_{k-j-1}(\alpha_{k-j})= \alpha_{j,k}   & &1\le j\le k\le\theta.
\end{align}
Let $| \, |:\Z^\theta\to\Z_2$ be the group map such that $|\alpha_i|=1$ if $i\in\{i_1,\dots,i_k\}$ and $|\alpha_i|=0$ otherwise, which defines the parity of the roots. $\cO(V)$ is the set of even roots, $\chi(\beta,\beta)=q^{\pm1}$ if $\beta\in\cO(V)$, and $\chi(\beta,\beta)=-1$ otherwise.

The algebra $U^+(V)$ is presented by generators $E_i$, $1\le i\le\theta$, and relations
\begin{align}
&E_iE_j=q_{ij} E_jE_i, &  &j-i \geq 2, \label{eq:relaciones U tipo A 1} \\
&[(\ad_c E_{i-1})(\ad_c E_i) E_{i+1}, E_i]_c, \quad E_i^2, & i &\in\{i_1,\dots,i_k\}, \label{eq:relaciones U tipo A 2} \\
&(\ad_c E_i)^2 E_{i\pm1}, & i & \notin\{i_1,\dots,i_k\}. \label{eq:relaciones U tipo A 3}
\end{align}
Let $E_{j,k} = (\ad_c E_j)\cdots (\ad_c E_{k-1})E_k= T_j\cdots T_{k-1}(E_k)$, $1\le j\le k\le\theta$. They are the generators of the PBW basis corresponding to the previous expression of $w$, ordered lexicographically:
$$ \alpha_{1,1}<\alpha_{1,2} <\cdots <\alpha_{1,\theta} <\alpha_{2,2} <\alpha_{2,3} <\cdots <\alpha_{\theta,\theta}. $$
By Remark \ref{rem:PRV=0 si alpha no cartan} $E_{j,k}^2=0$ for each odd root $\alpha_{j,k}$.

The computation of $\underline{\Delta}(E_{j,k})$ follows as for braidings of Cartan type with matrix $A_\theta$, see \cite[Lemma 6.5]{AS Pointed HA}.

\begin{lemma}\label{lema:coproducto Ejk, tipo A}
Let $1\le j\le k\le\theta$. Then
\begin{align}\label{eq:coproducto Ejk, tipo A}
\underline{\Delta}(E_{j,k})= E_{j,k}\ot 1+1\ot E_{j,k}+ \sum_{j\le \ell<k}(1-\widetilde{q_{l,l+1}}) E_{j,\ell} \ot E_{\ell+1,k}.
\end{align}
\end{lemma}
\pf
The proof is by induction on $k-j$, the case $k=j$ is trivial. Assume that \eqref{eq:coproducto Ejk, tipo A}
holds for $j+1$, $k$.
Notice that
\begin{align*}
\chi(\alpha_j,\alpha_{j+1,k})\chi(\alpha_{j+1,k},\alpha_j)&=\widetilde{q_{j,j+1}},&
\chi(\alpha_j,\alpha_{\ell+1,k})\chi(\alpha_{\ell+1,k},\alpha_j)&=1 \mbox{ if }\ell>j.
\end{align*}
Then we compute:
\begin{align*}
\underline{\Delta}(E_{j,k})& = \underline{\Delta}(E_j)\underline{\Delta}(E_{j+1,k})-\chi(\alpha_j,\alpha_{j+1,k}) \underline{\Delta}(E_{j+1,k})\underline{\Delta}(E_j) \\
&= E_{j,k}\ot 1+1\ot E_{j,k}+ (1-\widetilde{q_{j,j+1}}) E_j\ot E_{j+1,\ell} \\
& \qquad + \sum_{j+1\le \ell<k} (1-\widetilde{q_{l,l+1}}) \, E_{j,\ell} \ot E_{\ell+1,k}\\
& \qquad + \sum_{j+1\le \ell<k} (1-\widetilde{q_{l,l+1}})\chi(\alpha_j,\alpha_{j+1,k}) E_{j+1,\ell} \ot [E_j, E_{\ell+1,k}]_c,
\end{align*}
From \eqref{eq:relaciones U tipo A 1}, $[E_j, E_{\ell+1,k}]_c=0$ if $\ell>j$, which completes the proof.
\epf

We compute now $\underline{\Delta}(E_{j,k}^N)$ for each even root $\alpha_{j,k}$. Compare this with the case of Cartan braidings of type $A$ \cite[Lemma 6.9]{AS Pointed HA}.

\begin{proposition}
Let $1\le j\le k\le\theta$ be such that $\alpha_{j,k}\in\cO(V)$. Then
\begin{align}
\underline{\Delta}&(E_{j,k}^N)= E_{j,k}^N\ot 1+1\ot E_{j,k}^N \notag\\
& + \sum_{\ell: \, j\le \ell<k, \, \alpha_{j,\ell}\in\cO(V)} (1-\widetilde{q_{\ell,\ell+1}})^N \chi(\alpha_{j,\ell},\alpha_{\ell+1,k})^{\frac{N(N-1)}{2}} E_{j,\ell}^N \ot E_{\ell+1,k}^N.
\end{align}
\end{proposition}
\pf
By induction on $k-j$. The case $k=j$ is trivial. Assume that
\begin{align*}
\underline{\Delta}&(E_{j+1,k}^N)= E_{j+1,k}^N\ot 1+1\ot E_{j+1,k}^N \\
& + \sum_{\underset{\alpha_{j+1,\ell}\in\cO(\rho_j(V))} {\ell: \, j+1\le \ell<k,} } (1-\widetilde{\underline{q}_{\ell,\ell+1}})^N s_j^*\chi(\alpha_{j+1,\ell},\alpha_{\ell+1,k})^{\frac{N(N-1)}{2}} E_{j+1,\ell}^N \ot E_{\ell+1,k}^N.
\end{align*}
Then $ \Delta_j(E_{j+1,k}^N)=\underline{\Delta}(E_{j+1,k}^N)$ since $(\underline{\Delta}\ot\id)\underline{\Delta}(E_{j+1,k}^N)$ does not have terms $E_j^n$, $n\in\N$.
If $\ell>j$, then $\widetilde{\underline{q}_{\ell,\ell+1}}= \widetilde{q_{\ell,\ell+1}}$, and
\begin{itemize}
\item $\alpha_{j+1,\ell}\in\cO(\rho_j(V))$ if and only if $s_j(\alpha_{j+1,\ell})=\alpha_{j,\ell}\in\cO(V)$,
\item $s_j^*\chi(\alpha_{j+1,\ell},\alpha_{\ell+1,k})=\chi(s_j(\alpha_{j+1,\ell}),s_j(\alpha_{\ell+1,k}))=\chi(\alpha_{j,\ell},\alpha_{\ell+1,k})$.
\end{itemize}
We apply now Theorem \ref{thm:factorization delta Ti}:

\begin{align*}
\underline{\Delta}(E_{j,k}^N & )= \underline{\Delta}\circ T_j (E_{j+1,k}^N) = \Rf_i^{\rho_j(V)}\circ(T_j\ot T_j)\circ \Delta_j(E_{j+1,k}^N) \\
&=
\Rf_i^{\rho_j(V)} \Big( E_{j,k}^N\ot 1+1\ot E_{j,k}^N\Big) \\
& + \sum_{\underset{\alpha_{j,\ell}\in\cO(V)} {\ell: \, j+1\le \ell<k,} } (1-\widetilde{q_{\ell,\ell+1}})^N \chi(\alpha_{j,\ell},\alpha_{\ell+1,k})^{\frac{N(N-1)}{2}} \Rf_i^{\rho_j(V)} \Big(E_{j,\ell}^N \ot E_{\ell+1,k}^N \Big).
\end{align*}
$\Rf_i^{\rho_j(V)}\Big(E_{j,\ell}^N \ot E_{\ell+1,k}^N \Big)=\Big(E_{j,\ell}^N \ot E_{\ell+1,k}^N \Big)$ since $E_{\ell+1,k}^N \triangleleft\Ef_j^{(n)} =0$ for all $n\geq 1$, $\ell>j$. Also $\Rf_i^{\rho_j(V)} \left( E_{j,k}^N\ot 1 \right)=E_{j,k}^N\ot 1$. By Theorem \ref{thm:Z+ es subalg hopf U+} it suffices to compute $E_{j,k}^N \triangleleft\Ef_j^{(nN)}$. As $\underline{\Delta}$ is $\N_0$-graded, $E_{j,k}^N \triangleleft\Ef_j^{(nN)}=0$ if $n>1$, so we compute $X_N$, where $E_j^N\ot X_n$ is a summand in the expression of $\underline{\Delta}(E_{j,k}^N)=\underline{\Delta}(E_{j,k})^N$. By Lemma \ref{lema:coproducto Ejk, tipo A},
\begin{align*}
E_j^N\ot X_n & = (1-\widetilde{q_{j,j+1}})^N (E_j \ot E_{j+1,k})^N \\
& = (1-\widetilde{q_{j,j+1}})^N  \chi(\alpha_j,\alpha_{j+1,k})^{\frac{N(N-1)}{2}} \, E_j^N \ot E_{j+1,k}^N,
\end{align*}
so we have
$$ \Rf_i^{\rho_j(V)} \left( 1\ot E_{j,k}^N \right)= 1\ot E_{j,k}^N + (1-\widetilde{q_{j,j+1}})^N  \chi(\alpha_j,\alpha_{j+1,k})^{\frac{N(N-1)}{2}} \, E_j^N \ot E_{j+1,k}^N, $$
which completes the proof. \epf

\subsection{Braidings of type $\mathfrak{br}(2;5)$}
Now fix $\theta=2$, $\zeta$ a root of unity of order 5, $(q_{ij})_{i,j\in\I}$, $(r_{ij})_{i,j\in\I}$ two matrices such that
\begin{align*}
q_{11}&=\zeta, & \widetilde{q_{12}}&=\zeta^2, & q_{22}&=r_{22}=-1, & r_{11}&=-\zeta^3, & \widetilde{r_{12}}&=\zeta^3.
\end{align*}
Let $\chi,\psi:\Z^2\times\Z^2\to\ku^\times$ be the bicharacters defined by $(q_{ij})_{i,j\in\I}$, $(r_{ij})_{i,j\in\I}$, respectively, and $V$, $W$ the corresponding braided vector spaces.
They correspond to \cite[Table 1, row 13]{H-classif RS} and are related with the Lie superalgebra $\mathfrak{br}(2;5)$ over a field of characteristic 5  \cite{AA-WGCLSandNA}. Their generalized Cartan matrices are
$ C^{V}= \left( \begin{array}{cc} 2 & -3 \\ -1 & 2 \end{array}\right)$, $ C^W= \left( \begin{array}{cc} 2 & -4 \\ -1 & 2 \end{array}\right)$, and $\rho_1(V)=V$, $\rho_1(W)=W$, $\rho_2(V)=W$. The elements of maximal length are $(\sigma_1^{V}\sigma_2)^4$, $(\sigma_1^W\sigma_2)^4$, and the positive roots $\Delta_+^{V}$, $\Delta_+^W$ are, respectively,
\begin{align*}
\big\{ \alpha_1,3\alpha_1+\alpha_2,2\alpha_1+\alpha_2,5\alpha_1+3\alpha_2,3\alpha_1+2\alpha_2,4\alpha_1+3\alpha_2,\alpha_1+\alpha_2,
\alpha_2 \big\}, \\
\big\{\alpha_1,4\alpha_1+\alpha_2,3\alpha_1+\alpha_2,5\alpha_1+2\alpha_2,2\alpha_1+\alpha_2,3\alpha_1+2\alpha_2,\alpha_1+\alpha_2,
\alpha_2 \big\}.
\end{align*}
The sets of Cartan roots $\cO(V)$, $\cO(W)$ are, respectively,
\begin{align*}
&\left\{\alpha_1,2\alpha_1+\alpha_2,3\alpha_1+2\alpha_2,\alpha_1+\alpha_2\right\}, &
&\left\{\alpha_1,3\alpha_1+\alpha_2,2\alpha_1+\alpha_2,\alpha_1+\alpha_2 \right\}.
\end{align*}
The algebra $U^+(V)$ is presented by generators $E_1,E_2$ and relations
\begin{align*}
&E_2^2, &  &(\ad_cE_1)^4E_2, &  &[E_{1112},E_{112}]_c, & &\left[E_{4\alpha_1+3\alpha_2},E_{12}\right]_c.
\end{align*}
The algebra $U^+(W)$ is presented by generators $\Eb_1,\Eb_2$ and relations
\begin{align*}
&\Eb_2^2, &  &(\ad_c\Eb_1)^5\Eb_2, &  &\left[\Eb_1,\Eb_{3\alpha_1+2\alpha_2}\right]_c+r_{12}\Eb_{2\alpha_1+\alpha_2}^2, & &\left[\Eb_{3\alpha_1+2\alpha_2},\Eb_{12}\right]_c.
\end{align*}
Here $E_{i_1\cdots i_kj}=(\ad_c E_{i_1})\dots(\ad_c E_{i_{k}})E_j$, and analogous notation with $\Eb$.

\begin{lemma}\label{lema:root vectors br(2,5)}
The root vectors of $\cO(V)$, $\cO(W)$ are, respectively,
\begin{align*}
E_{\alpha_1}&=E_1, & E_{3\alpha_1+2\alpha_2}&= q_{12}^3\zeta^3(1-\zeta^3)^5(1+\zeta) [E_{112},E_{12}]_c, \\
E_{2\alpha_1+\alpha_2}&=q_{21}(1-\zeta^3)E_{112}, & E_{\alpha_1+\alpha_2}&=q_{21}^2\zeta(1-\zeta^3)^7(1+\zeta)^3 E_{12}, \\
\Eb_{\alpha_1}&=\Eb_1, & \Eb_{2\alpha_1+\alpha_2}&= -r_{12}^2\zeta^2(1+\zeta)^2(1-\zeta^3)^4 \Eb_{112}, \\
\Eb_{3\alpha_1+\alpha_2}&=r_{21}(1-\zeta^2)\Eb_{1112}, & \Eb_{\alpha_1+\alpha_2}&=r_{12}^2\zeta(1+\zeta)^3(1-\zeta^3)^{10} \Eb_{12}.
\end{align*}
\end{lemma}
\pf
First we compute
\begin{align*}
E_{2\alpha_1+\alpha_2}&=T_1^{V} T_2(\Eb_1)=T_1^{V}\left( E_{21} \right) = -q_{21}\zeta  \partial_1^L( E_{1112})=q_{21}(1-\zeta^3)E_{112}.
\end{align*}
A similar computation gives the expression of $\Eb_{3\alpha_1+\alpha_2}$. Using this,
\begin{align*}
E_{3\alpha_1+2\alpha_2} &=T_1^{V} T_2 (\Eb_{3\alpha_1+\alpha_2}) =(1-\zeta^2)^2 T_1^{V}([[E_{21},[E_{21},E_1]_c]_c) \\
&= [-L_1\cdot \partial_1^L(E_{112}), L_1\cdot \partial_1^L(L_1\cdot \partial_1^L(E_{112}))]_c\\
&= q_{12}^3\zeta^3(1-\zeta^3)^5(1+\zeta) [E_{112},E_{12}]_c.
\end{align*}
The proof for the remaining root vectors is analogous by iterating compositions of $T_1^{V} T_2$ and $T_1^W
 T_2$.
\epf

Now we compute the coproduct of powers of root vectors.

\begin{lemma}\label{lemma:coproducto br(2,5)}
$E_{\alpha_1}^5, E_{\alpha_1+\alpha_2}^{10}$, $\Eb_{\alpha_1}^{10}, \Eb_{\alpha_1+\alpha_2}^5$ are primitive, and
\begin{align*}
\underline{\Delta}(E_{2\alpha_1+\alpha_2}^{10})&= E_{2\alpha_1+\alpha_2}^{10}\ot 1+1\ot E_{2\alpha_1+\alpha_2}^{10} + \frac{q_{21}^{30}E_{\alpha_1}^{10}\ot E_{\alpha_1+\alpha_2}^{10}}{(1-\zeta^3)^{55}(1+\zeta)^{25}}   \\
& \qquad \qquad + \frac{q_{21}^{30}}{(1-\zeta^3)^{10}} E_{\alpha_1}^5 \ot E_{3\alpha_1+2\alpha_2}^5, \\
\underline{\Delta}(E_{3\alpha_1+2\alpha_2}^5)&= E_{3\alpha_1+2\alpha_2}^5  \ot 1+1\ot  E_{3\alpha_1+2\alpha_2}^5+ \frac{q_{12}^{15}E_{\alpha_1}^5 \ot E_{\alpha_1+\alpha_2}^{10}}{(1-\zeta^3)^{40}(1+\zeta)^{15}} ,\\
\underline{\Delta}(\Eb_{3\alpha_1+\alpha_2}^5)&= \Eb_{3\alpha_1+\alpha_2}^5  \ot 1+1\ot \Eb_{3\alpha_1+\alpha_2}^5 + \frac{r_{21}^{35}\zeta^2 \Eb_{\alpha_1}^{10} \ot \Eb_{\alpha_1+\alpha_2}^5}{(1-\zeta^3)^{40}(1+\zeta)^5}  , \\
\underline{\Delta}(\Eb_{2\alpha_1+\alpha_2}^{10})&=  \Eb_{2\alpha_1+\alpha_2}^{10} \ot 1+1\ot \Eb_{2\alpha_1+\alpha_2}^{10} - \frac{r_{21}^{45} \Eb_{\alpha_1}^{10} \ot \Eb_{\alpha_1+\alpha_2}^{10}}{(1-\zeta^3)^{50}}   \\
& \qquad \qquad - \frac{r_{12}^5(1+\zeta)^{10}}{(1-\zeta^3)^{10}} \Eb_{3\alpha_1+\alpha_2}^5 \ot  \Eb_{\alpha_1+\alpha_2}^5.
\end{align*}
\end{lemma}
\pf
We know that $E_{\alpha_1}^5$, $\Eb_{\alpha_1}^{10}$ are primitive. Now $E_{\alpha_1+\alpha_2}^{10}$, $\Eb_{\alpha_1+\alpha_2}^5$ are also primitive since $\underline{\Delta}$ is $\N_0^\theta$-graded and $Z^+(V)$ is a Hopf subalgebra by Theorem \ref{thm:Z+ es subalg hopf U+}. Also $\underline{\Delta}(E_{3\alpha_1+2\alpha_2}^5)$ is the sum of the two terms $\Eb_{3\alpha_1+\alpha_2}^5  \ot 1$, $1\ot \Eb_{3\alpha_1+\alpha_2}^5$ plus $c \,E_{\alpha_1}^5 \ot E_{\alpha_1+\alpha_2}^{10}$ for some $c\in \ku$ since the first term of the coproduct should be a power of previous roots by Proposition \ref{prop:lado izquierdo Delta Z(chi)} and then we use again Theorem \ref{thm:Z+ es subalg hopf U+}. As
$$  \underline{\Delta}([E_{112},E_{12}]_c^5)=(1\ot [E_{112},E_{12}]_c +E_1\ot \partial_1^L([E_{112},E_{12}]_c)+... )^5,$$
and $\partial_1^L([E_{112},E_{12}]_c)=-\zeta^3(1-\zeta^3)(1+\zeta)^2 E_{12}$, such $c$ is computed from $(E_1 \ot E_{12}^2)^5=q_{21}^{10} E_1^5\ot E_{12}^{20}$ using Lemma \ref{lema:root vectors br(2,5)}.
Similarly $\underline{\Delta}(E_{2\alpha_1+\alpha_2}^{10})$ can have two extra summands $a E_{\alpha_1}^{10}\ot E_{\alpha_1+\alpha_2}^{10}$ and $b E_{\alpha_1}^5 \ot E_{3\alpha_1+2\alpha_2}^5$. As $\partial_i^L(E_{112})=(1+\zeta)(1-\zeta^3)E_{12}$, the first appears from
$$ (E_1\ot E_{12})^{10} = q_{21}^{45} E_1^{10}\ot E_{12}  ,$$
while the second appears from the product
$$ ((1\ot E_{112})(E_1\ot E_{12}))^{5} = q_{21}^{15} E_1^5 \ot (E_{112}E_{12})^5 ,$$
The calculus of coproducts of the other two expressions is similar.
\epf

\begin{remark}\label{rem:accion br(2,5)}
The coproduct of the power root vectors in Lemma \ref{lemma:coproducto br(2,5)} seems close to the algebra of functions of $B_2$, and the Lie algebra $\g_0$ corresponding to the Lie superalgebra $\mathfrak{br}(2;5)$ in characteristic 5 is of type $B_2$. The PBW bases of $U(V)$ and $\u(V)$ resemble the PBW bases of the enveloping and the restricted enveloping algebras of the Lie superalgebra of $\mathfrak{br}(2;5)$ respectively, since the generators have the same degree and the corresponding height.
\end{remark}

\end{document}